\documentclass{amsart}

\usepackage{amsmath}
\usepackage{amsfonts}
\usepackage{amssymb}
\usepackage{latexsym}
\usepackage{graphicx}
\usepackage{enumerate}
\usepackage{multirow}
\usepackage{subfigure}

\newtheorem{theorem}{Theorem}
\newtheorem{proposition}{Proposition}
\newtheorem{assumption}{Assumption}
\newtheorem{lemma}[theorem]{Lemma}
\newtheorem{cor}[theorem]{Corollary}

\theoremstyle{definition}
\newtheorem{remark}{Remark}

\newcommand{\ep}{\epsilon}

\newcommand{\be}{\begin{equation}}
\newcommand{\ee}{\end{equation}}
\newcommand{\bes}{\begin{equation*}}
\newcommand{\ees}{\end{equation*}}
\newcommand{\R}{{\bf{R}}}

\newcommand{\ds}{\displaystyle}

\newcommand{\epi}{\| a^{-1} \|_{L^\infty_T}}

\newcommand{\abcd}{{a \tilde{b} \tilde{c} \tilde{d}}}
\newcommand{\epj}{{\epsilon_j}}
\newcommand{\epjp}{{\epsilon_{j'}}}

\begin{document}
\title[Dispersion vs. anti-diffusion]{Dispersion vs. anti-diffusion: well-posedness in 
variable coefficient and quasilinear  equations of  KdV-type.}
\author{David M. Ambrose}
\address{Department of Mathematics, Drexel University, 3141 Chestnut St., Philadelphia, PA 19104, USA}
\email{ambrose@math.drexel.edu}
\author{J. Douglas Wright}
\address{Department of Mathematics, Drexel University, 3141 Chestnut St., Philadelphia, PA 19104, USA}
\email{jdoug@math.drexel.edu}

\begin{abstract}
We study the well-posedness of the initial value problem on periodic intervals for linear and quasilinear evolution equations 
for which the leading-order terms
have three spatial derivatives.  In such equations, there is a competition between the dispersive effects which
stem from the leading-order term, and anti-diffusion which stems from the lower-order terms with two spatial derivatives.
We show that the dispersive effects can dominate the backwards diffusion:
we find a condition which guarantees well-posedness of the initial value problem for  
linear, variable coefficient equations of this kind, even when such anti-diffusion is present.
In fact, we show that even in the presence of localized backwards diffusion, the dispersion will
in some cases lead to an overall effect of parabolic smoothing. 
By contrast, we also show that when our condition is violated, the backwards diffusion can dominate the dispersive 
effects, leading to an ill-posed initial value problem.  We use these results on linear evolution equations as a guide when
proving well-posedness of the initial value problem for some quasilinear equations which also exhibit this competition
between dispersion and anti-diffusion: a Rosenau-Hyman compacton equation, the Harry Dym equation, and 
equations which arise in the numerical analysis of finite difference schemes for dispersive equations.  For these
quasilinear equations, the well-posedness theorem requires that the initial data be uniformly bounded away from zero.
\end{abstract}

\maketitle

\section{Introduction}
We study a family of evolution equations,
$$u_{t}=au_{xxx}+bu_{xx}+\phi,$$
where $a$ and $b$ are functions and $\phi$ represents a collection of lower-order terms, and
where $x \in X = [0,M]$ and we assume $u$ satisfies periodic boundary conditions.
We establish conditions on the coefficient functions $a$ and $b,$ as well as conditions on $\phi,$
which ensure that the initial value problem for this equation is well-posed in Sobolev spaces.  This is
a delicate question because if there exist any $x$ and $t$ for which $b<0,$ then formally, this term acts
(locally, at least) as an anti-diffusion or backwards heat operator.  The conditions we find on $a$ and $b$
indicate to what extent the dispersive effects of the leading-order term can compensate for the presence
of this anti-diffusion.  In a special case (in which the coefficient functions depend only on $x$), we also show that if our 
condition is violated, then the initial value problem is ill-posed.  We treat this problem for linear equations, for which
$a$ and $b$ are functions of $x$ and $t,$ and also certain quasilinear equations, for which $a$ and $b$ depend on $u$
and $u_{x}.$

Our original motivation for this work is to study the question of well-posedness or ill-posedness for the
the Rosenau-Hyman $K(2,2)$ equation,
\be\label{K22}
u_t = 2 uu_{xxx} + 6 u_x u_{xx} + 2 uu_x.
\ee
$K(2,2)$ is a degenerate version of the Korteweg-de Vries equation and famously
supports compactly supported traveling waves, dubbed compactons \cite{Rosenau-Hyman},
as well as other non-analytic traveling waves \cite{liOlverRosenau}.
Many papers in the literature make numerical simulations of solutions of the $K(2,2)$ initial value problem,
focusing primarily on the study of interactions of these compactons \cite{chertock}, \cite{deFrutos},
\cite{levyShuYan}, \cite{rus-pade}, \cite{rus-selfsimilar}.
While the existence theory
for this equation remains unsettled, there is strong evidence that 
$K(2,2)$ is ill-posed for 
 data which crosses zero \cite{ASWY}.
 Here we study a complementary
problem: Is the equation well-posed when the initial condition is bounded away from zero?
The answer, as we show below, is yes.

In order to show that $K(2,2)$ is well-posed, we first study an associated
linear problem:
\be\label{CKS}
u_t = a(x,t) u_{xxx} + b(x,t) u_{xx} + c(x,t) u_x + d(x,t) u + e(x,t)
\ee
where once again $x \in X = [0,M]$
and $u$ (and the coefficient functions) satisfy periodic boundary conditions.
We assume that the coefficient functions are defined for all $t \in [0,T],$
where $T > 0$.  

By way of analogy, supposing that $u(x,0)$ is bounded away from zero in $K(2,2)$ we should study
the situation where $a(x,t)$ is likewise bounded from zero.  On the other hand, $u_x(x,0)$ will change
signs for any non-trivial periodic initial data and so we should not impose strong conditions on $b(x,t)$.
In particular we must allow $b(x,t)$ to be negative, at least for some $x$.

Equations like \eqref{CKS}, and quasilinear variants, have been studied previously, though the focus has  been
on the case where $X = \R$ instead of the periodic problem we study here.
The earliest work \cite{Craig-Kappeler-Strauss} requires that $b \ge 0;$
that is to say,
there is no destabilizing ``backwards diffusion" effect from the term $b(x,t) u_{xx}$.
It is well-known that if $a(x,t)$ is identically zero and there is an interval on which $b < 0$, then these destabilizing
effect are so pathological that 
\eqref{CKS} is ill-posed  for data in Sobolev spaces.
{
In some sense, the results of this paper are an investigation into the extent to which 
{dispersive} effects from $a(x,t) u_{xxx}$
ameliorate the catastrophically destabilizing effects caused by a {backwards}
diffusive term.}

This question was answered for $X = \R$ in \cite{Akhunov}, and, in short,
the answer is that if $b(\cdot,t) \in L^1$ uniformly in $t$ and $a(x,t)$ is bounded away from zero,
then the dispersive effects completely dominate the backwards diffusion and
the equation (including quasilinear versions) is well-posed.  The issues for the periodic setting are
necessarily different; the argument cannot be based on integrability (the integrability would be considered over
the real line rather than a periodic interval, and a nontrivial periodic function is not integrable on the real line).
The answers that we find for the periodic case include a well-posedness theorem, a parabolic smoothing result, and 
an ill-posedness theorem for (\ref{CKS}), when certain conditions (primarily on $a$ and $b$) are satisfied.

While we were motivated by the well-posedness question for the $K(2,2)$ equation, 
our method of proof applies to other quasilinear evolution equations as well.  We do not
attempt to give the most general theorem that we can on well-posedness of quasilinear evolution
equations with nonlinear dispersion and anti-diffusion, but instead we consider some additional 
equations which have appeared in the literature.  One of these is the Harry Dym equation,
$
u_t = u^3 u_{xxx}.
$  Other examples come from the numerical analysis of finite difference schemes for certain dispersive
partial differential equations. 
We are able to prove well-posedness theorems for
such equations for uniformly positive Sobolev initial data.

The remainder of this paper is organized as follows.  In Section 2 we study the linear problem \eqref{CKS}
and provide sufficient conditions on the coefficients so that the equation is well-posed. The main tools here
are a change of dependent variables followed by energy estimates.  See Theorem \ref{linear main result}
below for the precise statement of this result.
In Section 3 we use the linear energy estimates as a inspiration for developing analogous estimates
for $K(2,2)$.  These estimates are sufficient for proving that $K(2,2)$ is well-posed.
See Corollary \ref{K22exist}
below for the precise statement of this result.  In Section 4, we briefly discuss additional quasilinear equations which
can be shown to be well-posed for positive initial data in Sobolev spaces.  
We mention that all of our well-posedness proofs largely follow the lines of Chapter 3 of \cite{majdaBertozzi}, in which
the Navier-Stokes equations are shown to be well-posed in Sobolev spaces.
In Section 5, we give a theorem on singularity formation:
if positive solutions of these quasilinear equations ever touch down and obtain a value of zero, then the solutions must blow up
in the space $H^{4}.$

\section{The linear problem}

\subsection{Main ideas and instructive examples}
The following heuristics and special cases  provide some 
crucial ideas for how we proceed. The first special case we consider is
\be\label{death}
u_t = a_0 u_{xxx} - b_0 u_{xx} \quad {\textrm{with}} \quad x \in [0,2 \pi]
\ee
where $a_0$ and $b_0$ are positive and constant and $u$ satisfies
periodic boundary conditions.
The equation can be solved {\it via} Fourier series:
$$
u(x,t) = \sum_{k \in {\bf Z}} 
 \exp(-i a_0 k^3 t) \exp(b_0 k^2 t) e^{i k x} \widehat{u}(k,0).
$$
Here $\widehat{u}(k,0)$ are the coefficients of the Fourier series
expansion of the initial data.
A routine calculation shows that the right hand side of this equation,
viewed as a map from $\R^+ \times H^s$ into $H^s$ does not
depend continuously on $u(x,0)$. This is due to the fact 
that 
$\exp(b_0 k^2t )$ grows extremely fast for large wave numbers $k$.
As such, \eqref{death} is ill-posed
in precisely the same way that backwards heat equation is.

If $b$ is uniformly negative, this example illustrates that 
there is no chance for the 
dispersive effects to play a role (above, we had $b(x,t)=-b_{0}<0$). 
However, if $b$ is localized and negative, there is
reason to think that the dispersion can arrest blow up,
as the following formal heuristic demonstrates.
Specifically, consider the following:
\be\label{life}
u_ t= a_0 u_{xxx} + b(x) u_{xx} \quad {\textrm{with}} \quad x \in [-M,M]
\ee
where 
$$
b(x) = \begin{cases}-b_0<0,& \text{ for $x \in [0,L]$}\\ 0, &\text{ otherwise}.\end{cases}
$$
Here $0 < L  \ll  M$.
 For simplicity, assume $a_0 < 0$.
Consider
a ``wave packet" of amplitude $\mu_0$ and wave number $k_0$
which is initially located to the left of $x = 0$. Such a wave packet will
move with group speed associated to the equation $u_t = a_0 u_{xxx}$,
which is to say with speed
$$c_g = 3|a_0| k_0^2.$$ 
When the packet enters the interval
$[0,L]$, say at $t=0$, it will be subjected to the effects of the backwards
heat term, which will cause the amplitude to grow like $$\mu_0 \exp(b_0 k_0^2 t),$$
as was the case in \eqref{death}.
The wave packet will leave the interval $[0,L]$ at time $t_{exit} = L/c_g$ and 
thus the final amplitude of the packet is 
\be\label{life2}
  \mu_0 \exp\left( b_0 k_0^2 t_{exit} \right)=
\mu_0 \exp\left(b_0 k_0^2 {L \over 3 |a_0| k_0^2}\right)=
  \mu_0 \exp\left({L \over 3} \left\vert  {b_0 \over a_0}\right \vert\right).
\ee
See Figure 1.
The growth
of the amplitude during the wave packet's passage through the destabilizing
zone is independent of its wave number---the extremely rapid growth of
high wave numbers is held in check by the fact that these highly oscillatory wave packets
move very fast.

\begin{figure}\label{heuristic}
\center{
\includegraphics[width=!]{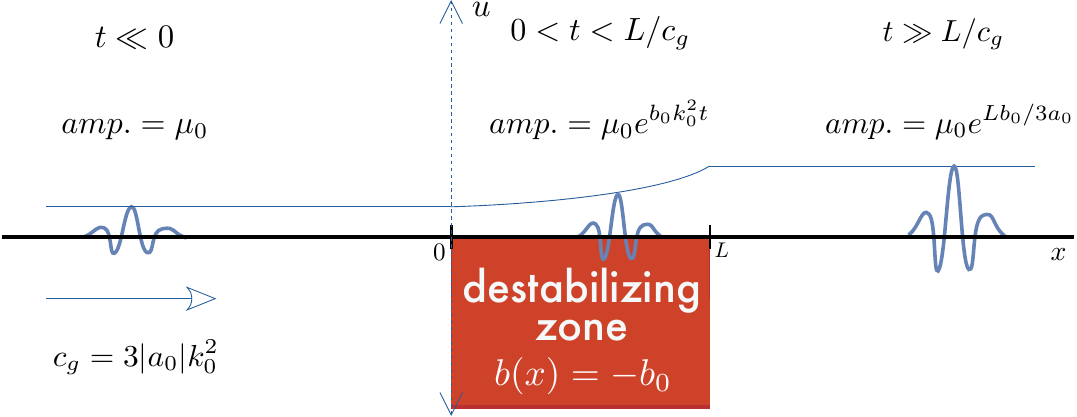} }
\caption{Dispersion {\it{vs}} localized backwards diffusion in \eqref{life}.}
\end{figure}

Therefore, if we presume that $b(x,t)$ is not everywhere negative, this heuristic
indicates a plausible mechanism by which dispersion can arrest blowup. 
Note that equation \eqref{life2} points to the importance of the ratio 
$$
\delta(x,t):=\ds {b(x,t) \over |a(x,t)|}
$$
in determining the growth of the amplitude. We call $\delta(x,t)$ the {\it modified
diffusion} for the equation, for reasons which will become clear below.

We make this formal argument---which is closely
related to the geometric optics heuristics described in \cite{Craig-Goodman}---precise in this article. Specifically, we give sufficient conditions
on the coefficient functions $a(x,t),\dots$, $e(x,t)$ so that \eqref{CKS} is well-posed.
In what follows $W^{k,p} = W^{k,p}(X)$ is the usual Sobolev space of periodic functions
on $X$, and $H^n=H^n(X) = W^{n,2}$.  Also, for any Banach space of functions defined
on $X$, we set
$$
B_T := L^\infty([0,T];B).
$$

We now state our assumptions on the coefficient functions:

\begin{assumption}\label{smoothness}
{\bf The coefficients are sufficiently regular/integrable:}
Fix $n \in {\bf N}$, with $n \ge 4$. Then
\begin{itemize}
\item $a(x,t) \in C([0,T];C^{n+3})$, 
\item $a_t(x,t) \in C([0,T];C^n)$,
\item $b(x,t) \in C([0,T];C^{n+2})$, 
\item $b_t(x,t) \in C([0,T];C^n)$,
\item $c(x,t) \in C([0,T];C^{n})$, 
\item $d(x,t) \in C([0,T];C^{n})$, and 
\item $e(x,t) \in L^\infty([0,T];H^n)$.
\end{itemize}
\end{assumption}

\begin{assumption}\label{uniform}{\bf The dispersion does not degenerate:}
There exists a constant $a_0$ with the property that  for all $x \in X$
and $t \in [0,T]$:
$$
0 < a_0 \le |a(x,t)|.
$$
\end{assumption}

%

\begin{assumption}
\label{modified}
{\bf The average modified diffusion is non-negative:}
$$
\bar{\delta}(t) : = {1 \over M}\int_0^{M} {b(y,t) \over |a(y,t)| }dy \ge 0
$$
for all $t \in [0,T]$.
\end{assumption}

\begin{remark}\label{what K means}
Let 
\begin{multline}\nonumber
\tilde{K}:=M+T+n+\|a\|_{W^{n+3,\infty}_T}+\|a_t\|_{W^{1,\infty}_T}+\epi \\
+\|b\|_{W^{n,\infty}_T}+\|b_t\|_{L^\infty_T}
+ \|c\|_{W^{n,\infty}_T}+ \|d\|_{W^{n,\infty}_T}+ \|e\|_{H^{n}_T}.
\end{multline}
Throughout this Section, any constant $K$ which appears in an
estimate has the property that $K>0$ and $K$ depends only on $\tilde{K}$. For instance,
one can show that
$$
\bar{\delta}(t) + |{d \over dt} \bar{\delta}(t)| \le K.
$$
\end{remark}

Our first main result is:

\begin{theorem}\label{linear main result}
Fix $n\in{\bf{N}}$ such that $n \ge 4$ and take 
Assumptions \ref{smoothness}, \ref{uniform}  and \ref{modified} to be true.
Then there is a continuous map
$$
\Psi: H^n \longrightarrow L^\infty([0,T];H^n)
$$
such that for all $u_{0}\in H^{n},$
$$
\| \Psi(u_0) \|_{ L^\infty([0,T];H^n) } \le K \left( \| u_0 \|_{H^n} + 1\right),
$$
and such that the following are true for $u = \Psi(u_0)$:
\begin{enumerate}
\item For all $0\le n'<n$,  $u \in\ds  C([0,T];H^{n'})$ and $\|u(t)\|_{H^{n'}} \le K \left( \| u_0 \|_{H^n} +1 \right)$,
\item $u(x,0) = u_0$, 
\item $u(x,t)$ solves \eqref{CKS} for $t \in [0,T]$ and
\item $u(x,t)$ is the only function for which (2) and (3) hold simultaneously.
\end{enumerate}
\end{theorem}
We have two extensions of this if we have more information about $\bar{\delta}(t)$.
First, 
\begin{theorem}\label{parabolic smoothing}
If, in addition to Assumptions \ref{smoothness}, \ref{uniform}  and \ref{modified}, we have for all $t \in [0,T]$
$$
\bar{\delta}(t) \ge \delta_0 > 0
$$ 
then the solution map $\Psi$ from Theorem \ref{linear main result} enjoys a parabolic smoothing property. That is, for all $u_0 \in H^n$,
$$
\| \Psi(u_0) \|_{ L^2([0,T];H^{n+1}) } \le K \| u_0 \|_{H^n}.
$$ Moreover, if Assumption \ref{smoothness} holds with $n$ replaced by $n+1$, then
$u=\Psi(u_0) \in C([0,T]; H^{n})\cap C((0,T];H^{n'})$ for all $n' \in [0,n+1)$.
\end{theorem}

\begin{theorem}\label{reversible}
If, in addition to Assumptions \ref{smoothness}, \ref{uniform}  and \ref{modified}, we have for all $t \in [0,T]$
$$
\bar{\delta}(t) = 0, 
$$ 
then then the solution map $\Psi$ from Theorem \ref{linear main result} satisfies
$u=\Psi(u_0) \in C([0,T];H^n)$.
\end{theorem}

\begin{remark}
To demonstrate how surprising our results are, consider the following examples, all posed for $x \in [0,2 \pi]$.
Clearly, the initial value problem for the equation
$$
u_t = \sin(x) u_{xx}
$$
is ill-posed. 
However, if we add even a very small amount of dispersion, for instance, 
$$
u_t = {1 \over 10^{23}} u_{xxx} + \sin(x) u_{xx}
$$
then the initial value problem for the equation
is well-posed since the average modified diffusion in this case is equal to zero (since $\sin(x)$ has average value zero).

Similarly, for any $\varepsilon$ such that $0<\varepsilon<1,$ the initial value problem for the equation
$$
u_{t} = \left(\varepsilon + \sin(x)\right)u_{xx}
$$
is again ill-posed.  However, if we consider the equation with some dispersion, say
$$ u_{t} = \frac{1}{10^{23}}u_{xxx} + \left(\varepsilon + \sin(x)\right)u_{xx},$$
then since
$
\bar{\delta} = 10^{23} \varepsilon > 0,
$
the initial value problem is well-posed.
Morever, Theorem \ref{parabolic smoothing} states that even for $a(x,t)$
 ludicrously
small, as in this case, so long as $\varepsilon > 0$, the equation is essentially parabolic even
though, on a sizable fraction of the interval, there is a backwards heat operator acting.
\end{remark}

\begin{remark}
One can ask which of the assumptions are necessary as well as sufficient. 
We do not have definitive answers at this time,
but some relevant information is known.
The  results in \cite{ASWY} and \cite{Craig-Goodman}
strongly indicate that if $a(x,t)$ crosses zero for some $x$
with non-zero derivative, then \eqref{CKS} will be ill-posed.
The regularity/integrability conditions in Assumption 
\ref{smoothness} can almost certainly be weakened, but 
doing so is likely technical and not terribly interesting
at this time. 
 A quick calculation shows
that sending $t \to -t$ in \eqref{CKS} results in an equation of the same
form as \eqref{CKS} but reverses the sign of $\bar{\delta}$.
The parabolic smoothing property in Theorem \ref{parabolic smoothing}
indicates that it is highly unlikely that \eqref{CKS} would be well-posed
backwards in time when $\bar{\delta}(t) > 0$. These two facts together
indicate that Assumption \ref{modified} is almost certainly necessary.
We confirm this in the special case when the  coefficients functions are time independent 
in Theorem \ref{ill posed} below. 
%
\end{remark}

\begin{remark}\label{a is positive}
Assumption \ref{uniform} allows us, without loss of generality, to take
$$
a(x,t) \ge a_0 > 0.
$$
And so, from this point on, $a(x,t)$ is presumed to be positive.
\end{remark}

\subsection{The change of variables and key estimate} Our goal is to find an {\it a priori} estimate
for solutions $u(x,t)$ of \eqref{CKS} in the space $H^n$.  It turns out that it
this more easily done after changing dependent variables.
Toward this end, we let $g_n(x,t)$
be function which has the following three properties:
\begin{enumerate}
\item [(C1)] There is a constant $k_g > 1$ such that for all $x$ and $t,$
$$
k_g^{-1} \le g_n(x,t) \le k_g,
$$
\item[(C2)] $g_n \in C([0,T];C^{n+3})$ and $\partial_t g_n \in C([0,T];C^{n})$, and
\item[(C3)] $g_n$ and its first $n$ derivatives satisfy periodic boundary conditions.
\end{enumerate}
Then define $v(x,t)$ {\it via}
\be\label{w def}
u(x,t) := g_n(x,t)v(x,t).
\ee
An intelligent choice for $g_n$ will allow us to find 
an energy estimate for $v,$ and we will make such a choice at the appropriate place in our proof.
 The following Lemma,
whose proof is obvious, implies
we will have a comparable estimate for $u$.
\begin{lemma}\label{u or v} Assume $g_n$ satisfies (C1)-(C3). Then
$u \in H^n_T$ if and only if $v \in H^n_T$.
That is, there is a constant $C_g \ge 1$, depending only on $\| g_{n} \|_{W^{n,\infty}_T}$ and $\| g_{n}^{-1} \|_{L^\infty_T}$, such that
$$
C_g^{-1} \| u \|_{H^n_T} \le \| v\|_{H^n_T} \le C_g \| v \|_{H^n_T}.
$$
\end{lemma}

Substituting \eqref{w def} into \eqref{CKS} yields the following after a routine
calculation:
\begin{equation}\label{CKS2}
v_t =  a v_{xxx} + \tilde{b} v_{xx} + \tilde{c} v_x + \tilde{d}v + \tilde{e}
\end{equation}
where
\be\label{COV}
\begin{split}
\tilde{b}(x,t) &= 3 a {\partial_x g_n\over g_n}  +b \\
\tilde{c}(x,t) &= 3a {\partial_x^2 g_n \over g_n}  + 2b {\partial_x g_n\over g_n}+c\\
\tilde{d}(x,t) & = -{\partial_t g_n\over g_n} + a {\partial_x^3 g_n \over g_n} + b {\partial_x^2 g_n \over g_n} 
+ c{\partial_x g_n \over g_n} + d\\
\tilde{e}(x,t) &= {e \over g_n}.
\end{split}
\ee

\begin{remark}
The regularity assumptions in (C2) are made so that $\tilde{b}$, $\tilde{c}$, $\tilde{d}$ and $\tilde{e}$ 
are in the same spaces as $b$, $c$, $d$ and $e$.
\end{remark}

Note that the right hand side of \eqref{CKS2} is of the same form as that of \eqref{CKS}.
Motivated by this, define
the operator 
$$
L_{ABCD}:= A \partial_x^3 + B \partial_x^2 + C \partial_x + D 
$$
where $A$, $B$, $C$ and $D$ are functions  which satisfy periodic boundary conditions in $x$ .
The following Lemma contains the heart of our energy estimate:
\begin{lemma}\label{heart lemma} Fix $n \in {\bf N}$ with $n \ge 4$. Suppose that $A,B,C,D \in W^{n,\infty}$.
Then there exists $\kappa>0$ (which depends only on $n$, $M$ and  the $W^{n,\infty}$ norms of $A,B,C$ and $D$) such that, for all $\phi \in H^{n}$:
\be\label{heart1}
\int_X (\partial_x^n \phi) \partial_x^n \left( L_{ABCD} \phi \right)\ dx
-\int_X \left[ \left( {3 \over2} -n \right)A_x - B\right] [\partial_x^{n+1}\phi]^2\ dx \le 
 \kappa \| \phi \|^2_{H^n}.
\ee
\end{lemma}
\noindent{\bf Proof:}
We prove the result for smooth functions $\phi$ only.  The case for  $\phi \in H^{n}$  follows by a typical density argument.
From the defintion of  $L_{ABCD}$, we have:
\begin{equation}\label{start}\begin{split}
\int_X (\partial_x^n \phi) \partial_x^n \left( L_{ABCD} \phi \right)\ dx&=
\int_X (\partial_x^n \phi)  \partial_x^n \left( A \phi_{xxx} + B \phi_{xx} + C\phi_{x} + D\phi\right)\ dx\\&=
I_A + I_B + I_C + I_D,
\end{split}
\ee
where  $I_A$, $I_B$, $I_C$ and $I_D$ are defined in the obvious way. We estimate each 
in turn.

First, the product rule gives:
\be\nonumber
\begin{split}
I_A := &\int_X (\partial_x^n \phi) \partial_x^n(A \phi_{xxx})\ dx \\=&
 \int_X (\partial_x^n \phi) A \partial_x^{n+3} \phi\ dx
 +\int_X n (\partial_x^n \phi) A_x \partial_x^{n+2} \phi\ dx
 \\&+ \int_X {n \choose 2} (\partial_x^n \phi) A_{xx} \partial_x^{n+1} \phi\ dx
 + \sum_{j = 0}^{n-3}\int_X {n \choose j} (\partial_x^n \phi) (\partial_x^{n-j} A) \partial_x^{j+3} \phi\ dx.
 \end{split}
 \ee
We call the four terms on the right hand side $I_A^1$, $I_A^2$, $I_A^3$, $I_A^4$ respectively.

$I^4_A$ is controlled with the Cauchy-Schwarz inequality
and other elementary tools, since at most $n$ derivatives of $\phi$ appear:
 \be\label{standard trick 1}
I_A^4 :=  \sum_{j = 0}^{n-3}\int_X {n \choose j} (\partial_x^n \phi)(\partial_x^{n-j} A) \partial_x^{j+3} \phi\ dx
 \le \kappa \| \phi \|^2_{H^n}.
 \ee
$I_A^3$  is handled using the observation that $\partial_x^n \phi \partial_x^{n+1} \phi = \ds {1 \over 2} \partial_x [\partial_x^n \phi]^2$. Thus an integration by parts gives:
\begin{multline}\label{standard trick 2}
 I_A^3 : = \int_X {n \choose 2} (\partial_x^n \phi) A_{xx} \partial_x^{n+1} \phi\ dx
= {1 \over 2} {n \choose 2} \int_X A_{xx} \partial_x [\partial_x^n \phi]^2\ dx
\\=-{1 \over 2} {n \choose 2} \int_X A_{xxx} [\partial_x^n \phi]^2\ dx \le \kappa \| \phi \|^2_{H^n}.
\end{multline}

A quick calculation shows that, for any smooth function $\psi$,
$$
\ds \psi \psi_{xxx} = {1 \over 2} \partial_x^3 [\psi^2] - {3 \over 2} \partial_x [\psi_x]^2.
$$
Applying this with  $\psi = \partial_x^n \phi$ in $I^1_A$
gives:
$$
I_A^1:=\int_X A (\partial_x^n \phi)  \partial_x^{n+3} \phi\ dx=
{1 \over 2} \int_X A \partial_x^3 [ \partial_x^n \phi]^2\ dx - {3\over2} \int_X A \partial_x [\partial_x^{n+1} \phi ]^2\ dx.
$$
Integrating by parts three times in the first term and just once in the second term gives:
\begin{equation}\nonumber
I_A^1 = -{1 \over 2} \int_X A_{xxx}  [ \partial_x^n \phi]^2\ dx + {3\over2} \int_X A_x [\partial_x^{n+1} \phi ]^2\ dx
 \le  {3\over2} \int_X A_x [\partial_x^{n+1} \phi ]^2\ dx + \kappa \| \phi\|^2_{H^n}.
\end{equation}

Now we consider
$$\ds
I^2_A:=\int_X n A_x (\partial_x^n \phi) \partial_x^{n+2} \phi\ dx.
$$
For smooth functions $\psi$, we have
\be\label{trick}
\psi \psi_{xx} = {1 \over 2} \partial_x^2[\psi^2] - [\psi_x]^2.
\ee
If we use this with $\psi = \partial_x^n \phi$, then
$$
 I_A^2 
 = {1 \over 2} \int_X n A_x \partial_x^2 [ \partial_x^n \phi]^2\ dx - \int_X n A_x [\partial_x^{n+1}\phi]^2\ dx.
$$
Integrating by parts twice in the first term gives:
$$
I_A^2 
 = {1 \over 2} \int_X n A_{xxx} [ \partial_x^n \phi]^2 dx - \int_X n A_x [\partial_x^{n+1}\phi]^2\ dx 
 \le - \int_X n A_x [\partial_x^{n+1}\phi]^2\ dx  + \kappa \| \phi \|^2_{H^n}.
$$

Adding the estimates for $I_A^{1}$, $I_A^{2}$, $I_A^3$ and $I_A^4$ yields:
\be\label{IA}
I_A \le \left( {3\over2} -n \right) \int_X A_x [\partial_x^{n+1} \phi ]^2\ dx +  \kappa \| \phi \|^2_{H^n}.
\ee

Now we estimate $I_B$. As with $I_A$, the product rule gives:
\begin{equation*}
\begin{split}
I_B : = &\int _X (\partial_x^n \phi)  \partial_x^n \left( B \phi_{xx} \right)\ dx\\
= & \int_X (\partial_x^n \phi) B \partial_x^{n+2} \phi\ dx
 +\int_X n (\partial_x^n \phi) B_x \partial_x^{n+1} \phi\ dx
\\& + \sum_{j = 0}^{n-2}\int_X {n \choose j} (\partial_x^n \phi)(\partial_x^{n-j} B) \partial_x^{j+2} \phi\ dx.
\end{split}
\end{equation*}
Using the same steps as in \eqref{standard trick 1} and \eqref{standard trick 2}, we have
$$
\int_X n (\partial_x^n \phi) B_x \partial_x^{n+1} \phi\ dx\\
 + \sum_{j = 0}^{n-2}\int_X {n \choose j} (\partial_x^n \phi)(\partial_x^{n-j} B) \partial_x^{j+2} \phi\ dx
\le \kappa \| \phi \|^2_{H^n}.
$$
Notice that the other term in $I_B$, namely $\ds \int_X B (\partial_x^n \phi)  \partial_x^{n+2} \phi\ dx$,
is the same as $I^2_A$ but with $B$ replacing $n A_x$. So using \eqref{trick}
and following the steps used to estimate $I_A^2$ we get:
$$
\ds \int_X B (\partial_x^n \phi)  \partial_x^{n+2} \phi dx
  \le - \int_X  B [\partial_x^{n+1}\phi]^2 dx  + \kappa \| \phi \|^2_{H^n}.
$$
Combining these estimates, we have
\be \label{IB}
I_B \le  - \int_X B [\partial_x^{n+1}\phi]^2\ dx  + \kappa \| \phi \|^2_{H^n}.
\ee
Controlling $I_C$ and $I_D$ amounts to applying
the processes in \eqref{standard trick 1} and \eqref{standard trick 2} appropriately.
One finds
$$
I_C + I_D \le \kappa \| \phi \|^2_{H^n}.
$$
Adding this to \eqref{IA} and \eqref{IB} then rearranging terms in \eqref{start} completes the proof of the lemma. \hfill$\blacksquare$

\subsection{Selection of $g_n$}
Notice that if
$
\ds \left( {3 \over2} -n \right)A_x - B \le 0,
$
that the second term on the left hand side of \eqref{heart1} is non-negative, so dropping
it gives
\be\label{heart}
\int_X (\partial_x^n \phi) \partial_x^n \left( L_{ABCD} \phi \right)\ dx \le 
 \kappa \| \phi \|^2_{H^n}.
\ee
This fact 
will lead to an {\it a priori} $H^n$ estimate for solutions $v$ of $\eqref{CKS2}$
provided $$\ds \left( {3 \over2} -n \right)a_x - \tilde{b} \le 0.$$ Given the definition of $\tilde{b}$ this is the same
as requiring that the following differential inequality is satisfied:
\be\label{C4}
 \left( {3 \over 2} -n \right) a_x - b - 3 a {\partial_x g_n \over g_n}\le 0 \quad \text{for all $x \in X$ and $t \in [0,T]$}.
\ee
At first glance, the 
 most straightforward way to meet \eqref{C4} is to choose $g_n$
so that that $\ds \left( {3 \over 2} -n \right) a_x - b - 3 a {\partial_x g_n \over g_n} = 0$.
That is to say, take $g_n$ to be
$$
\left(a(x,t)\right)^{1/2 - n/3} \exp\left[ -{1 \over 3} \int_0^x {b(y,t) \over a(y,t)}\ dy \right].
$$
The problem with doing this is the fact that this function does not meet periodic boundary conditions,
as
$\ds
 \int_0^M {b(y,t) \over a(y,t)}\ dy
$
is not necessarily zero. 
Since our estimates require multiple integrations by parts, it is absolutely critical that $g_n$ be periodic.
Assumption \ref{modified} is made precisely so that we may choose $g_n$ to be periodic while simultaneously
satisfying \eqref{C4}. (It is this fact that $g_n$ must be periodic
which makes the conditions on $a$ and $b$ more restrictive than in the case when one studies \eqref{CKS}
for $x \in \R$.)

To find a function $g_{n}$ which satisfies our conditions, we begin by taking $\bar{\delta}(t)$ as in Assumption \ref{modified}. 
Then take $g_n$
to be a non-zero solution of 
\be\label{g eqn}
{\partial_x g \over g} = -{1 \over 3} \left( {b \over a} - \bar{\delta} \right) + \left({1 \over 2} - {n \over 3} \right) {a_x \over a}.
\ee
Specifically, define
\be \label{gn}
g_n(x,t) := \left(a(x,t)\right)^{1/2 - n/3} \exp\left[ -{1 \over 3} \int_0^x\left(  {b(y,t) \over a(y,t)} - \bar{\delta}(t) \right)\ dy \right].
\ee
Then we have the following lemma:
\begin{lemma}\label{choice}
Take Assumptions \ref{smoothness}, \ref{uniform} and \ref{modified} as given
and define $g_n$ as in \eqref{gn}. Then $g_n$ satisfies (C1), (C2), (C3) and \eqref{C4}.
\end{lemma}
\noindent{\bf Proof:}
Assumptions \ref{smoothness} and \ref{uniform} imply $a$ is bounded away from zero and also bounded above,
so the same is true for $\left(a(x,t)\right)^{1/2 - n/3}$. Likewise $b$ is bounded above, as is $\bar{\delta}$, and so
the exponential is bounded above and away from zero. Thus we have (C1) for an appropriate constant $k_g$.
The regularity of $g_n$ can be ascertained from \eqref{gn}. Essentially $g_n$ has the same regularity as $a$ (due
to the prefactor of $a^{1/2-n/3}$) and one derivative smoother than $b$ (since $b$ appears in \eqref{gn} only as 
$\ds\int_0^x \frac{b(y)}{a(y)}\ dy$).
Thus Assumption \ref{smoothness}
gives (C2). Observe that, for any $t \in [0,T]$ the definition of $\bar{\delta}$ and the fact that $a$ is periodic gives:
\begin{multline*}
g_n(M,t) = \left(a(M,t)\right)^{1/2 - n/3} \exp\left[ -{1 \over 3} \int_0^M\left(  {b(y,t) \over a(y,t)} - \bar{\delta}(t) \right)\ dy \right]\\
=\left(a(M,t)\right)^{1/2 - n/3} \exp\left[ -{1 \over 3}\left( M \bar{\delta}(t) - M\bar{\delta}(t) \right) \right]\\
=\left(a(M,t)\right)^{1/2 - n/3} = \left(a(0,t)\right)^{1/2 - n/3} =g_{n}(0,t).
\end{multline*}
Thus $g_n$ is periodic. The same reasoning shows that its derivatives are periodic and we have (C3).
Finally, rearranging the terms in \eqref{g eqn} gives
\be\label{choice 2}
 \left( {3 \over 2} -n \right) a_x - b - 3 a {\partial_x g_n \over g_n} = - a \bar{\delta}(t).
\ee
Since $a > 0$ and since Assumption \ref{modified} tells us that $\bar{\delta} \ge 0$, we have \eqref{C4}.
\hfill$\blacksquare$

\subsection{Well-posedness of \eqref{CKS2}}
Now that we have selected $g_n$ we can prove that \eqref{CKS2}, and consequently \eqref{CKS}, is well-posed.
We will first prove uniqueness and continuous dependence of solutions, and we will then prove existence.
For all of these, a uniform bound for solutions will be helpful.
Lemma \ref{choice} leads to the following energy estimate:
\begin{proposition}\label{energy}
Take Assumptions \ref{smoothness}, \ref{uniform} and \ref{modified} as given
and define $g_n$ as in \eqref{gn}.
If $v \in L^\infty([0,T];H^n)$ solves \eqref{CKS2} then
$$
{d \over dt} \| v (t) \|_{H^n}^2 \le K \| v (t) \|^2_{H^n}  + 2 (v(t) ,\tilde{e}(t))_{H^n}.
$$
\end{proposition}
\noindent{\bf Proof:} Following the time-honored tradition, we compute
$$
{d \over dt} \| v(t) \|^2_{H^n} =
{d \over dt} \| v(t) \|^2_{L^2}  + {d \over dt} \|\partial_{x}^{n}v(t) \|^2_{L^2} =
2\int_X v v_t\ dx + 2\int_X (\partial_x^n v) \partial_x^n v_t\ dx.
$$
Using the definition of $L_{ABCD}$, $\eqref{CKS2}$ can be rewritten as $v_t = L_{\abcd }v + \tilde{e}$.
Thus:
\begin{multline*}
{d \over dt} \| v(t) \|^2_{H^n} = 2\int_X v L_{\abcd} v\ dx +2 \int_X (\partial_x^n v) \partial_x^n \left(L_{\abcd} v \right)\ dx\\
+ 2\int_X v \tilde{e}\ dx +2 \int_X (\partial_x^n v) \partial_x^n \tilde{e}\ dx\\
=  2\int_X v L_{\abcd} v\ dx +2 \int_X (\partial_x^n v) \partial_x^n \left(L_{\abcd} v \right)\ dx + 2 (v(t) ,\tilde{e}(t))_{H^n}.
\end{multline*}

Since $n \ge 4$ and $L_\abcd$ takes three derivatives:
$$
2\int_X v L_{\abcd} v\ dx \le K \|v(t)\|_{H^n}^2. 
$$
Our choice of $g_n$ and Lemma \ref{choice} allow us to use
\eqref{heart}, thus:
$$
2 \int_X \partial_x^n v \partial_x^n \left(L_{\abcd} v \right)\ dx \le  K \| v(t) \|^2_{H^n}.
$$
\hfill$\blacksquare$

This proposition immediately implies: 
\begin{cor} \label{cdoic}
(Uniqueness and continuous dependence on initial conditions)
Take Assumptions \ref{smoothness}, \ref{uniform} and \ref{modified} as given
and define $g_n$ as in \eqref{gn}. 
Suppose that $v_1$, $v_2 \in H^n_T$ are solutions of \eqref{CKS2}. Then
for all $t \in [0,T]$,
$$
\| v_1(t) - v_2(t) \|_{H^n} \le K \| v_1(0) - v_2(0) \|_{H^n}.
$$
In particular, if $v_1(0)  = v_2(0)$ then $v_1 \equiv v_2$.
\end{cor}
\noindent{\bf Proof:}
Let $\Delta = v_1 - v_2$. Clearly 
$$
\Delta _t = L_{\abcd} \Delta
$$
and $\Delta(0) = v_1(0) - v_2(0)$. Apply Proposition \ref{energy} with $\tilde{e} = 0$ to get
$\ds {d \over dt} \| \Delta (t) \|^2_{H^n} \le K  \| \Delta(t) \|_{H^n}^2$. Gronwall's inequality then gives
the  estimate.
\hfill$\blacksquare$

All that remains is to show that \eqref{CKS2} does indeed have solutions. Here is our result:
\begin{proposition}\label{v wp}Take Assumptions \ref{smoothness}, \ref{uniform} and \ref{modified} as given
and define $g_n$ as in \eqref{gn}. 
For all $v_0 \in H^n$, with $n \ge 4$, there exists a function $v^*$ with the following properties:
\begin{itemize}
\item $v^* \in L^\infty([0,T];H^n)$ with $\| v \|_{H^n_T} \le K \left(\|v_0\|_{H^n} + 1\right)$,
\item for all $0 \le n' < n$, $v^* \in C([0,T];H^{n'})$ with $\| v \|_{H^{n'}_T} \le K \left(\|v_0\|_{H^n} + 1\right)$,
\item $v^*(x,0) = v_0$ and
\item $v^*(x,t)$ solves \eqref{CKS2}.
\end{itemize}
\end{proposition}

\begin{remark}
This proposition, Corollary \ref{cdoic}, and Lemma \ref{u or v} imply the main result, Theorem \ref{linear main result}.
\end{remark}

\noindent{\bf Proof:} (Proposition \ref{v wp})
We first regularize \eqref{CKS2}. Let
${\mathcal{J}_\ep} \psi := \eta_\epsilon * \psi$ where $*$ is the periodic convolution
and $\eta_\epsilon$ is  ``Evans' standard mollifier" \cite{evans}.  For utility we record the following
properties of ${\mathcal{J}_\ep}$:
\begin{theorem}\label{mollifier stuff}
\begin{enumerate}
\item For all $\ep > 0$, $s_1, s_2 \in \R$, $\| {\mathcal{J}_\ep} \|_{H^{s_1} \to H^{s_2} } \le C(\ep,s_1,s_2)$.
\item If $\phi \in H^s$, then $\lim_{\epsilon \to 0^+} \| {\mathcal{J}_\ep} \phi - \phi \|_{H^s} = 0$.
\item For all $\ep > 0$, $s \in \R$, $\| {\mathcal{J}_\ep} \|_{H^s \to H^s} \le C(s)$. Specifically, $C(s)$
does not depend on $\epsilon$.
\item For all $\ep > 0$, $s \in \R$, and $\phi,\psi \in H^s$, $({\mathcal{J}_\ep} \phi,\psi)_{H^s} = ( \phi,{\mathcal{J}_\ep}\psi)_{H^s}$.
\item ${\mathcal{J}_\ep} \partial_x = \partial_x {\mathcal{J}_\ep}$.
\end{enumerate}
\end{theorem}
Consider now the following regularized initial vale problem:
\be
\label{CKSmol}
v^\ep_t = {\mathcal{J}_\ep}\left(  a {\mathcal{J}_\ep}  v^\ep_{xxx} + \tilde{b} {\mathcal{J}_\ep} v^\ep_{xx} + \tilde{c} {\mathcal{J}_\ep} v^\ep_x + \tilde{d} \mathcal{J}_\ep v^\ep + \tilde{e} \right)
\ee
with
\be\label{ic}
v^\ep(x,0) = v_0 \in H^n.
\ee
Due to the presence of the  mollifiers, the 
 right hand side of \eqref{CKSmol} defines a Lipschitz
 map from $H^n$ to itself.  Since ${\mathcal{J}_\ep}$ commutes with derivatives, we can rewrite the right hand side of \eqref{CKSmol} as:
$$
{\mathcal{J}_\ep} L_\abcd {\mathcal{J}_\ep} v^\epsilon+ {\mathcal{J}_\ep} \tilde{e}.
$$ 
So we have from the standard existence theory 
 for ODEs on a Banach spaces:
 \begin{lemma} \label{CKSmol exist}
 For all $\epsilon > 0$ there exist $T_\ep > 0$ and
$$
 v^\ep(x,t) \in C^1([0,T_\ep];H^n)
 $$
which solves \eqref{CKSmol} and \eqref{ic}.  
\end{lemma}
Here, $T_\ep > 0$ depends, possibly badly, on $\epsilon$. 
We now prove an analog of Proposition \ref{energy} for solutions of \eqref{CKSmol}. Specifically,
we have, just as in the proof of Proposition \ref{energy}:
\begin{multline}\nonumber
{d \over dt} \| v^\epsilon(t) \|^2_{H^n}=  2\int_X v^\epsilon {\mathcal{J}_\ep} \left( L_{\abcd}  {\mathcal{J}_\ep} v^\ep\right)\ dx 
\\
+2 \int_X \partial_x^n v^\ep \partial_x^n {\mathcal{J}_\ep} \left(L_{\abcd} {\mathcal{J}_\ep} v^\ep \right)\ dx + 2 (v^\ep(t) ,
{\mathcal{J}_\ep} \tilde{e}(t))_{H^n}.
\end{multline}
Using the fact that ${\mathcal{J}_\ep}$ is self-adjoint (Theorem \ref{mollifier stuff}, part (4)) gives:
\begin{multline}\label{pickup}
{d \over dt} \| v^\epsilon(t) \|^2_{H^n}=  2\int_X {\mathcal{J}_\ep} v^\ep  \left( L_{\abcd}  {\mathcal{J}_\ep} v^\ep\right)\ dx 
\\
+2 \int_X \partial_x^n {\mathcal{J}_\ep} v^\ep \partial_x^n  \left(L_{\abcd} {\mathcal{J}_\ep} v^\ep \right)\ dx 
+ 2 (\mathcal{J}_\ep v^\ep(t) ,\tilde{e}(t))_{H^n}.
\end{multline}
The first term is easily bounded by $K \| {\mathcal{J}_\ep} v^\epsilon \|^2_{H^n}$, since $n \ge 4$. We can apply Lemma 
\ref{heart lemma} and the
estimate \eqref{heart}  to the second term, with $\phi = {\mathcal{J}_\ep} v^\ep$.  (The mollifiers
were arranged precisely in \eqref{CKSmol} so that we could do this.) We then have:
$$
{d \over dt} \| v^\epsilon(t) \|^2_{H^n} \le K \| {\mathcal{J}_\ep} v^\epsilon(t) \|^2_{H^n}+2 (\mathcal{J}_\ep v^\ep(t) ,\tilde{e}(t))_{H^n}.
$$
Note that the constant $K$ here (and any such constants which appear below) do not depend on $\epsilon.$

Since $\| {\mathcal{J}_\ep} \phi \|_{H^n} \le C(n) \| \phi \|_{H^n}$ where $C(n)$ does not depend on $\epsilon$, we have, with the Cauchy-Schwarz inequality:
$$
{d \over dt} \| v^\epsilon(t) \|^2_{H^n} \le K \left( \| v^\epsilon(t) \|^2_{H^n}+ \| v^\ep(t)\|_{H^n} \right)\le K \left(1+ \| v^\epsilon(t) \|^2_{H^n}\right).
$$
Applying Gronwall's inequality to this gives, for all $t \in [0,T_\ep]$.
\be\label{gron}
\| v^\ep(t) \|^2_{H^n} \le (1 + \| v_0 \|_{H^n}^2) e^{Kt} - 1 \le K( \| v_0 \|_{H^n}^2  + 1).
\ee
Thus we have $\| v^\ep(T_\ep) \|_{H^n} < \infty$. This means that we can continue the solution; we can do so indefinitely 
given the form of \eqref{gron}.
  We conclude
that, for all $\epsilon$, 
\be\label{picard}
v^\ep \in C^1([0,T];H^n).
\ee
Which is to say that the solutions $v^\ep$ all exist for a common time interval.

Since $n \ge 4$, we have $H^n\subset  C^3$. Thus the estimate \eqref{gron} implies
that 
$$
\sup_{t \in [0,T]} \| v^\ep(t) \|_{C^3} \le K( \| v_0 \|_{H^n}  + 1).
$$
And since $v^\ep$ solves \eqref{CKSmol}, we have
$$
\sup_{t \in [0,T]} \|v^\ep_t(t) \|_{C} \le K\sup_{t \in [0,T]}  \| v(t) \|_{C^3} \le K( \| v_0 \|_{H^n}  + 1).
$$
In particular, there is a constant $K$ independent of $\epsilon$:
$$
\sup_{t \in [0,T]} \left( \| v^\ep_x \|_{C}  + \|v^\ep_t \|_{C}  \right) \le K( \| v_0 \|_{H^n}  + 1).
$$
This implies that $\left\{ v^\ep \right\}_{\ep > 0}$ is a uniformly bounded and equicontinuous
family of functions defined for $(x,t) \in [0,M] \times [0,T]$. $[0,M] \times [0,T] \in \R^2 $ is compact,
and so we apply the Arzela-Ascoli theorem to conclude that there is a function $v^* \in C([0,M] \times [0,T])$
and subsequence $\left\{\ep_j\right\}$ (with $\lim_{j \to \infty} \ep_j = 0$) such that
$$
\lim_{j \to \infty} \| v^{\ep_j} - v^*\|_{C([0,M] \times [0,T])}
=\lim_{j \to \infty} \| v^{\ep_j} - v^*\|_{L^\infty_T }=0.
$$
We have from interpolation, for all $n' \in [0,n)$:
\be\label{interp}
\| v^{\epj}(t) - v^{\epjp}(t) \|_{H^{n'}} \le \|v^{\epj}(t) - v^{\epjp}(t)\|^{q}_{L^2} \| v^\epj(t) - v^{\epjp}(t) \|^{r}_{H^n}
\ee
where $q = 1 - n'/n$ and $r = n'/n$. The estimate \eqref{gron} implies that:
$$
\| v^{\epj}(t) - v^{\epjp}(t) \|_{H^{n'}} \le K ( \| v_0 \|_{H^n}  + 1)^r \|v^{\epj}(t) - v^{\epjp}(t)\|^{q}_{L^2}.
$$
Since $[0,M]$ is compact,  $L^\infty \subset L^2$:
\be\label{cauchy}
\| v^{\epj} - v^{\epjp} \|_{H^{n'}_T} \le K ( \| v_0 \|_{H^n}  + 1)^r \|v^{\epj} - v^{\epjp}\|^{q}_{L^\infty_T}.
\ee
Of course, since $v^\epj \to v^*$ uniformly, $\left\{ v^\epj \right\}$ is a Cauchy
sequence in the $C([0,T];L^\infty)$ topology. 
Moreover, we know that for all $\epsilon$, $v^\epsilon \in C^1([0,T];H^{n'})$ by \eqref{picard}.
Thus \eqref{cauchy} implies $\left\{ v^\epj \right\}$ is a Cauchy sequence
in the $C([0,T];H^{n'})$ topology, and therefore convergent. 
Of course, the limit must be $v^*$.
And so, for all $n \in [0,n')$:
$$
v^* \in C([0,T];H^{n'}).
$$

Moreover, since the $v^\ep$ are in $C([0,T];H^n)$, 
for each $t \in [0,T]$ we have $\left\{ v^\epj(t) \right\}$ a bounded sequence in $H^n$. 
Thus this sequence has a weak limit in $H^n$. Call the limit $v^{**}(t)$. We have
$
v^{**} \in L^\infty([0,T];H^n).
$
An elementary argument using the fact that limits are unique shows that
$v^{**} = v^*$. But this then means
$$
v^* \in L^\infty([0,T];H^{n}).
$$
We claim that $v^*$ solves $\eqref{CKS2}$ with initial condition \eqref{ic}.
For all $\epsilon$, we have
\be\label{time int}
v^\epsilon(x,t) = v_0(x) + \int_0^t \left( {\mathcal{J}_\ep} L_\abcd {\mathcal{J}_\ep} v^\epsilon(x,s)+ {\mathcal{J}_\ep} \tilde{e}(s)\right) ds.
\ee
Recall that $L_\abcd$ takes three derivatives and $\ds\lim_{\ep \to 0^+} \| {\mathcal{J}_\ep} \phi - \phi \|_{H^s} =0 $
for any $\phi \in H^s$. Moreover
$v^\epj \to v^*$ as $j \to \infty$ in $C([0,T];H^{n'})$
for any $n' \in [0,n)$ with $n \ge 4$.
And so we have no difficulty passing to the limit in \eqref{time int} along the subsequence $\left\{ \epsilon_j \right\}$ to get
$$
v^*(x,t) = v_0(x) + \int_0^t \left( L_\abcd  v^*(x,s)+  \tilde{e}(s)\right) ds.
$$
This implies that $v^*$ sastifies \eqref{CKS2} and \eqref{ic}.
\hfill$\blacksquare$

We now prove the extensions.

\noindent{\bf Proof:} (Theorem \ref{parabolic smoothing}) Suppose that
\be\label{parab}
\bar{\delta}(t) \ge \delta_0 > 0.
\ee
Follow the proof of Theorem \ref{v wp} until \eqref{pickup}. Then $\eqref{gron}$
implies
$$
{d \over dt} \| v^\epsilon(t) \|^2_{H^n} \le  2 \int_X \partial_x^n {\mathcal{J}_\ep} v^\ep \partial_x^n  \left(L_{\abcd} {\mathcal{J}_\ep} v^\ep \right) dx +  K \left(\| v_0\|_{H^n}^2 + 1\right)
$$
which, upon subtracting from both sides, gives:
\begin{multline*}
{d \over dt} \| v^\epsilon(t) \|^2_{H^n} 
-2 \int_X \left[ \left( {3 \over2} -n \right)a_x - \tilde{b}\right] [\partial_x^{n+1}{\mathcal{J}_\ep} v^\ep]^2 \\
\le  2 \int_X \partial_x^n {\mathcal{J}_\ep} v^\ep \partial_x^n  \left(L_{\abcd} {\mathcal{J}_\ep} v^\ep \right) dx 
\\
-2\int_X \left[ \left( {3 \over2} -n \right)a_x - \tilde{b}\right] [\partial_x^{n+1}{\mathcal{J}_\ep} v^\ep]^2+  K \left(\| v_0\|_{H^n}^2 + 1\right).\end{multline*}
Applying \eqref{heart1} on the right hand side we get:
$$
{d \over dt} \| v^\epsilon(t) \|^2_{H^n} 
-2 \int_X \left[ \left( {3 \over2} -n \right)a_x - \tilde{b}\right] [\partial_x^{n+1}{\mathcal{J}_\ep} v^\ep]^2  \le K \left(\| v_0\|_{H^n}^2 + 1\right).
$$
Our choice \eqref{gn} for $g_n$ and \eqref{choice 2} give:
$$
{d \over dt} \| v^\epsilon(t) \|^2_{H^n} +
 2 \bar{\delta} (t)  \int_X a [\partial_x^{n+1}{\mathcal{J}_\ep} v^\ep]^2  \le  K \left(\| v_0\|_{H^n}^2 + 1\right). 
$$
Integrating both sides with respect to time gives us the following estimate:
$$
 \| v^\epsilon(t) \|^2_{H^n} +
 2 \int_0^t \bar{\delta} (s)  \int_X a(s) [\partial_x^{n+1}{\mathcal{J}_\ep} v^\ep(s)]^2 ds \le K \left( \|v_0\|^2_{H^n} + 1\right).
$$
Thus we have, for all $t \in [0,T]$:
$$
 \int_0^t \bar{\delta} (s)  \int_X a(s) [\partial_x^{n+1}{\mathcal{J}_\ep} v^\ep(s)]^2 ds \le K  \left( \|v_0\|^2_{H^n} + 1\right).
$$
Assumption \ref{uniform} and \eqref{parab} imply that
$$
\int_0^t \| \partial_x^{n+1} {\mathcal{J}_\ep} v^\ep(s)) \|_{L^2}^2 ds \le {1 \over \delta_0 a_0}  \int_0^t \bar{\delta} (s)  \int_X a(s) [\partial_x^{n+1}{\mathcal{J}_\ep} v^\ep(s))]^2 ds.
$$
Thus, for any $\epsilon > 0$:
$$
\int_0^T \| \partial_x^{n+1} {\mathcal{J}_\ep} v^\ep(s)) \|_{L^2}^2 ds \le  K  \left( \|v_0\|^2_{H^n} + 1\right).
$$
And so $\{ {\mathcal{J}_\ep} v^\ep \}$ is a bounded set in $L^2([0,T];H^{n+1})$. It then has a weakly convergent
subsequence, whose limit
we call $v^{***} \in  L^2([0,T];H^{n+1})$. Given that $v_\ep \to v^*$ strongly in $C([0,T];H^{n'})$ when $0 \le n' < n$,
it is straightforward to conclude that ${\mathcal{J}_\ep} v^\ep \to v^*$  strongly in $L^2([0,T];H^{n'})$. Uniqueness of limits thus
implies $v^* = v^{***}$. Thus $v^* \in  L^2([0,T];H^{n+1})$. As norms are lower semi-continuous with respect to weak limits we have
$$
\int_0^T \| v^*(s)\|^2_{H^{n+1}} ds \le \lim_{\ep \to 0}  \int_0^T \| {\mathcal{J}_\ep} v^\ep(s)) \|_{H^{n+1}}^2 ds \le  K  \left( \|v_0\|^2_{H^n} + 1\right).
$$

All that remains is to show that $v^* \in C([0,T];H^n)$.
Since $v^* \in L^2([0,T];H^{n+1})$, we see that for a.e. $t \in [0,T]$, $v^*(t) \in H^{n+1}$.
Let $\tau > 0$ be such a time. Then we can rerun the existence argument
in its entirety using $v^*(\tau) \in H^{n+1}$ as the initial condition, but
in the smoother space (hence the additional regularity requirements
in the statement of Theorem \ref{parabolic smoothing}).  Call this new solution $v^\star$. By Proposition \ref{v wp},
we have $v^\star \in C([\tau,T];H^{n'+1})$ for all $n' \in [0,n)$. By uniqueness,
it must be the case that $v^\star = v^*$ for $t \in [\tau,T]$. Lastly, we can conclude that
$$
\lim_{t \to 0^+} \| v^* - v_0\|_{H^n} = 0;
$$
the details of this step are the same as in the proof of Theorem \ref{K22 exist} below.
Therefore, we get $v^{*}\in C([0,T]; H^{n}),$ and since $\tau$ can be taken arbitrarily close to $0$, we also have
$v^* \in C((0,T];H^{n'+1})$ as claimed.
\hfill$\blacksquare$

\noindent{\bf Proof:} (Theorem \ref{reversible}) Suppose that
\be\label{ref}
\bar{\delta}(t)  = \frac{1}{M} \int_0^M {{b}(y,t) \over |{a}(y,t)|}\ dy= 0.
\ee
The crux of the argument here is that in this case  \eqref{CKS} is,
for all intents and purposes, reversible. Specifically, suppose that $u(x,t)$ solves
\eqref{CKS} with initial conditions $u_0(x)$. Then let $w(x,t) = u(M-x,t_0-t)$
where $t_0 \in (0,T]$. 
Then
\be\label{CKSrev}
w_t = \bar{a} w_{xxx} - \bar{b} w_{xx} + \bar{c} w_x - \bar{d} w + \bar{e}
\ee
where $\bar{a}(x,t) = a(M-x,t_0-t)$, $\bar{b}(x,t) = b(M-x,t_0-t)$ and so on. 
Solving \eqref{CKS} for $t$ backwards from $t_0$ is therefore equivalent to solving \eqref{CKSrev}
forward in time. 
The claim is that  \eqref{CKSrev} meets all the Assumptions, and is thus solvable.
Clearly the coefficients
on the right hand side of \eqref{CKSrev} satisfy the regularity requirements of Assumption \ref{smoothness}.
Likewise it is clear that $\inf_{x \in X} | \bar{a}(x,t) | = \inf_{x \in X} | a(x,t) |$ and so Assumption \ref{uniform} is also met.
Lastly,
$$
{1 \over M} \int_0^M {\bar{b}(y,t) \over |\bar{a}(y,t)|}\ dy=
{1 \over M} \int_0^M {{b}(M-y,t) \over |{a}(M-y,t)|}\ dy = {1 \over M} \int_0^M {{b}(y,t) \over |{a}(y,t)|}\ dy=0.
$$
Thus Assumption \ref{modified} is met, and so we apply Theorem \ref{linear main result} to get solutions of \eqref{CKSrev}.
In combination with the energy estimates we have, the reversibility implies 
$$
u \in C([0,T];H^{n});
$$
we omit the argument at present, as it is essentially the same as the argument in the proof of Theorem \ref{K22 exist} below.
\hfill$\blacksquare$

\subsection{Ill-posedness when $\bar{\delta} < 0$}

Let us define $L$ as:
$$
L u : = a(x) u_{xxx} + b(x) u_{xx} + c(x) u_x + d(x) u.
$$
We have:
\begin{theorem}\label{ill posed}
Suppose\footnote{We take smooth periodic coefficients
here for simplicity.} that $a$, $b$, $c$ and $d$ are in $C^\infty(X)$.
If 
$$
\inf_{x \in X} a(x) = a_0 > 0
\quad \text{and} \quad
\bar{\delta} ={1 \over M} \int_0^M {b(y) \over a(y)}\ dy \ne 0
$$
then there is a sequence of eigenvalues $\left\{ \lambda_j \right\}_{j \in {\bf N}}$of $L$ for which
$$
\lim_{j \to \infty}  -\bar{\delta}\   \Re(\lambda_j)= +\infty.
$$
\end{theorem}
Notice that if $\bar{\delta} < 0$ then this proposition implies that $L$ has eigenvalues with arbitrarily
large real part. And thus $u _ t = L u$ is ill-posed.

\noindent{\bf Proof:}
We treat the special case 
 when
 $$
 Lu = u_{xxx} + \bar{\delta}u_{xx} + c(x)u_x + d(x) u.
 $$
 The general case follows by a change of variables, which we postpone
 for the moment.
 This portion
of the proof is similar to the proof that uniformly elliptic operators
have eigenvalues which tend to infinity in \cite{evans}. Note that 
we are viewing $L$ as an unbounded linear operator from $L^2$ to itself
whose domain is $H^3$.

First we claim that $L$ has a sequence of eigenvalues $\left\{ \lambda_j \right\}_{j \in \bf{N}}$ for which
$$
\lim_{j \to \infty} | \lambda_j | = \infty.
$$
Since we are working on the periodic interval, $L$ is a relatively compact perturbation 
of $\partial_x^3$.  Thus we have $\sigma_{ess}(L) = \sigma_{ess}(\partial_x^3) = \left\{ \right\}$.
This implies that the resolvent set of $L$ is nonempty. Let $\omega$ be a point in the resolvent.
Let $K = (L-\omega)^{-1}$. We have $K$ is bounded linear map from $L^2$ to $H^3$.
The Rellich-Kondrachov theorem implies that $H^3$ is compactly embedded in $L^2$,
and thus we conclude $K$ is a compact map from $L^2$ to itself. This implies
that $K$ has a sequence of non-zero eigenvalues $\left\{\mu_j\right\}_{j \in {\bf N}}$ for which $\lim_{j \to \infty} \mu_j = 0$.
Since $K$ is bounded from $L^2$ to $H^3$ we have functions $u_j \in H^3$ so that $K u_j = \mu_j u_j$.
Since $K = (L-\omega)^{-1}$ a short calculation shows
$$
L u_j = \left(\omega + {1 \over \mu_j} \right) u_j.
$$
Thus $\lambda_j =\ds  \omega +  {1 \over \mu_j}$ is an eigenvalue of $L$.  Since the $\mu_j$
tend to zero, the $|\lambda_j|$ tend to infinity. Thus the claim is shown.

Next,
let $u$ and $\lambda$ be an eigenfunction/eigenvalue pair for $L$, with $\| u \|_{L^2} = 1$.
Since $Lu = \lambda u$ we have
\be \label{lambda1}
\lambda = \int_X Lu(x) \bar{u}(x) dx   = \int_X (u_{xxx} + \bar{\delta} u_{xx} + c(x) u_x + d(x) u ) \bar{u}\ dx.
\ee
Now, by integration by parts,
 $\ds \int_X u_{xxx} \bar{u}\ dx = - \int_X u \bar{u}_{xxx}\ dx.$ This implies that $\Re\ds \int_X u_{xxx} \bar{u}\ dx = 0$.
Likewise, integration by parts shows $\int_X \bar{\delta} u_{xx} \bar{u}\ dx = -\bar{\delta} \| u_x \|^2_{L^2}$.
And so taking the real part of \eqref{lambda1} gives:
\be\label{real part}
\Re(\lambda) = -\bar{\delta}\|u_x\|^2_{L^2} + \Gamma
\ee
where
$$
\Gamma :=  \Re\left[ \int_X\left(c(x) u_x + d(x) u \right) \bar{u}\ dx \right].
$$
A straightforward calculation shows
$$
\left| \Gamma \right| \le K\|u\|_{L^2}^2 = K
$$
where $K>0$ depends only on $c$ and $d$.
Notice that \eqref{real part} implies, when $\bar{\delta} < 0$ that
$
\Re(\lambda) \ge -K
$
for some constant $K$ which depends only $c$ and $d$. On the other hand
if $\bar{\delta} >0$ then we have $\Re(\lambda) \le K$. Which is to say
\be\label{real part 2}
\bar{\delta} \Re(\lambda) \le |\bar{\delta}| K.
\ee
for all eigenvalues $\lambda$.
Solving \eqref{real part} for $\|u_x\|_{L^2}$ gives:
\be\label{ux}
\|u_x\|_{L^2}^2 =-\bar{\delta}^{-1} \left(\Re(\lambda) - \Gamma\right) \le k_1 |\Re(\lambda)| + k_2
\ee
for positive constants $k_1$ and $k_2$ which do not depend on $\lambda$.

Taking imaginary parts of \eqref{lambda1} gives:
$$
\Im(\lambda) = \int_X u_{xxx} \bar{u}\ dx + \Omega
$$
where
$$
\Omega = \Im\left[ \int_X\left(c(x) u_x + d(x) u \right) \bar{u}\ dx \right].
$$
One can show that using the Cauchy-Schwarz inequality that:
$$
\Omega \le K \|u_x\|_{L^2} \|u\|_{L^2} = K \|u_x\|_{L^2}.
$$
where $K$ depends only on $c$. 

An integration by parts and Cauchy-Schwarz give 
$$\int_X u_{xxx} \bar{u}\ dx = -\int_X u_{xx} \bar{u}_x dx \le \|u_{xx}\|_{L^2} \|u_x\|_{L^2}.$$ So we have
$$
|\Im(\lambda)| \le \left( \|u_{xx}\|_{L^2} + K\right) \|u_x\|_{L^2}.
$$
Using \eqref{ux} gives:
\be\label{halfway done}
|\Im(\lambda)| \le \left( \|u_{xx}\|_{L^2} + K\right) \sqrt{k_1| \Re(\lambda)| + k_2}.
\ee

Now, since $Lu = \lambda u$ we have
$$
\lambda \int_X u \bar{u}_{xx}\ dx = \int_X \left( u_{xxx} + \bar{\delta} u_{xx} + c(x) u_x + d(x) u \right) \bar{u}_{xx}\ dx.
$$
An integration by parts shows that the left hand side equals $-\lambda \|u_x\|_{L^2}^2$. Therefore we have,
after a rearrangement of terms:
$$
\bar{\delta} \|u_{xx}\|_{L^2}^2 = - \lambda \|u_x\|^2_{L^2} - \int_X \left( u_{xxx}  + c(x) u_x + d(x) u \right) \bar{u}_{xx}\ dx.
$$
Integration by parts shows that $\ds \Re \int_X u_{xxx} \bar{u}_{xx}\ dx = 0$. Taking the real part of the last equation
then gives
$$
\bar{\delta} \|u_{xx}\|_{L^2}^2 = - \Re(\lambda) \|u_x\|^2_{L^2} +\Theta
$$
with 
$$
\Theta:=
- \Re\left[\int_X \left(   c(x) u_x + d(x) u \right) \bar{u}_{xx}\ dx\right].
$$
Integration by parts shows that
$$
\left| \Theta \right| \le K\left( \| u_x\|^2_{L^2}  + 1\right)
$$
where $K>0$ depends on $c$ and $d$. 
$$
\|u_{xx}\|_{L^2}^2 = -\bar{\delta}^{-1} \Re(\lambda) \| u_x\|^2_{L^2} + \bar{\delta}^{-1} \Theta
$$
and using the estimate on $\Theta$ and also \eqref{ux} we get:
\be
\|u_{xx}\|_{L^2}^2 \le k_3 \Re(\lambda)^2 + k_4
\ee
for positive constants $k_3$ and $k_4$ which do not depend on $\lambda$.
Inserting this into \eqref{halfway done} gives
\be\label{almost there}
|\Im(\lambda)| \le \left(\sqrt{  k_3 \Re(\lambda)^2 + k_4} +K\right)  \sqrt{k_1| \Re(\lambda)| + k_2}.
\ee
Now suppose that $\bar{\delta} < 0$. Then \eqref{real part 2} tells us that 
$
\Re(\lambda) \ge -K < 0.
$
This, with \eqref{almost there} implies that all the eigenvalues of $L$ are in the set
$$
\Sigma:=\left\{ z = x+iy \in {\bf C} : x > -K {\text{ and }} |y| \le  \left(\sqrt{  k_3 x^2 + k_4} +K\right)  \sqrt{k_1|x| + k_2}.\right\}
$$
But we also know that there is a sequence of eigenvalues of $L$ for which $|\lambda| \to \infty$. A sequence
which diverges to infinity with $\Sigma$ must have
$\Re(\lambda) \to + \infty$. The situation when $\bar{\delta} > 0$ is similar.

Thus we have shown our result for the case when $a = 1$ and $b(x) = \bar{\delta}$. A
change of variables much like the one we used in the study of \eqref{CKS} will allow us to reduce the general
case to this one, though there is one wrinkle.  First set
$$
\xi(x) = \int_0^x a^{-1/3}(y)\ dy. 
$$
so that $d \xi /dx = a^{-1/3}(x)$. By assumption that $a \ge a_0$ we see that $\xi$ is a strictly increasing function,
and thus invertible.
Then define $\tilde{v}(\eta)$ by $u(x) = \tilde{v}(\xi(x))$. A tedious but routine calculation shows that
$$
Lu = a(x) \left(   \xi'  \right)^3 \tilde{v}_{\eta \eta \eta} 
+ \left( 3 a(x) \xi' \xi'' + b(x) (\xi')^2 \right) \tilde{v}_{\eta \eta} + c_1(\eta) \tilde{v}_{\eta} + d_1(\eta) \tilde{v}.
$$
Formulae for the functions $c_1$ and $d_1$ are omitted, as they are complicated and not useful.

Our choice for $\xi$ was made so that $a(x) \left(   \xi'  \right)^3 = 1$. Likewise
$$ 
\left( 3 a(x) \xi' \xi'' + b(x) (\xi')^2 \right) = a^{-2/3}(x)\left(b(x) - a'(x) \right)=: B(\eta).
$$
Observe that
$$
\int_0^{\tilde{M}} B(\eta) d\eta = \int_0^M B(\xi(x)) \xi'(x) dx. 
$$
Here $\tilde{M} = \xi^{-1}(M)$. Since $\xi$ is increasing, $\tilde{M} > 0$.
And thus
\begin{multline}\nonumber
\int_0^{\tilde{M}} B(\eta) d\eta= \int_0^M a^{-2/3}(x)\left(b(x) - a'(x) \right)a^{-1/3} (x)dx \\
=  \int_0^M\left({ b(x) \over a(x)} - {a'(x) \over a(x) }\right)dx =\ds  M \bar{\delta}.
\end{multline}
For this reason, set $\ds \bar{B} = {M \over \tilde{M}} \bar{\delta}$. $\bar{B}$ has the same sign
as $\bar{\delta}$.
Define
$$
\tilde{L} \tilde{v} := \tilde{v}_{\eta \eta \eta} + B(\eta) \tilde{v}_{\eta \eta} + c_1(\eta) \tilde{v}_{\eta} + d_1(\eta) \tilde{v}.
$$
Since $\tilde{L}$ is obtained from $L$ by a change of variables, clearly they have the same eigenvalues.
Following \eqref{gn}, let $\tilde{v} = g v$ where 
$$
g = \exp\left( - {1 \over 3} \int_0^x\left( B(s) - \bar{B} \right) ds \right).
$$
A quick calculation shows that if $\tilde{L} \tilde{v} = \lambda \tilde{v}$ then
$$
v_{\eta \eta \eta} + {M \over \tilde{M}} \bar{\delta} v_{\eta \eta} + c_2(\eta) v_{\eta} + d_2(\eta) v = \lambda v
$$
where $c_2$ and $d_2$ are functions whose formulae we omit. The operator on the left hand side
is exactly of the form we first studied in the proof. Since $\ds{M \over \tilde{M}} \bar{\delta}$ has the same
sign as  $\bar{\delta}$, we are done.
\hfill$\blacksquare$

\section{Well-posedness of the $K(2,2)$ IVP}

Recall that the $K(2,2)$  equation \eqref{K22} is:
$$
u_{t}=2uu_{xxx}+6u_{x}u_{xx}+2uu_{x}.
$$
In this section we prove that this is well-posed for 
for uniformly positive $H^{n}$ initial data on the periodic interval $X$.
The method mirrors that which worked for \eqref{CKS}. Specifically,
we change variables and then prove energy estimates.  For the linear equation \eqref{CKS}
the change of variables was linear, and was given by \eqref{w def} and \eqref{gn}. 
It turns out
that applying \eqref{w def} and \eqref{gn} to a solution $u$ of \eqref{K22} does
in fact lead to a new formulation for \eqref{K22} in which one can prove
{\it a priori} estimates on the $H^n$ norm. This is done as follows.
Formally,
\eqref{K22} is an equation of the form \eqref{CKS} with $a = 2u$, $b = 6 u_x$ and $c = 2u$.
Inserting these choices for $a$ and $b$ into \eqref{gn} yields:
\be \label{gn2}
g_n(x,t) := \left(2u(x,t)\right)^{1/2 - n/3} \exp\left[ -{1 \over 3} \int_0^x\left(  {6u_x(y,t) \over 2u(y,t)} - \bar{\delta}(t) \right)\ dy \right]
\ee
where $\ds \bar{\delta}(t) = {1 \over M} \int_0^M  {6u_x(y,t) \over 2u(y,t)}\ dy$, as in Assumption \ref{modified}.
Now, $u_x/u = \partial_x(\ln(u))$ and we are assuming $u$ is periodic and bounded away from zero. 
Thus we have $\ds \bar{\delta}(t) =0$. With this, we can evaluate $g_n$ above to get:
$$
g_n(x,t) = C(t) u^{-1/2-n/3}
$$
for a function $C(t)$ which does not depend on $x$. For simplicity we take $C(t) = 1$.
Following \eqref{w def}, we set
$w=u/g_n$. All together, this suggests that changing variables in \eqref{K22} {\it via}
$$
w(x,t) := [u(x,t)]^{{9 + 2n \over 6}}
$$
will lead to estimates on the $H^n$ norm of $w$.  And so we set
\be\label{beta}
\beta := { 6\over9 + 2n }
\ee
and define $w$ as
$$
u = w^\beta.
$$
Note that for $n  \ge 0$, $\beta \in (0,1)$. This nonlinear change of variables is
completely modified in proving existence of solutions. 
It is worth remarking that while the linear estimates provided a roadmap
for finding this change of variables, our nonlinear results
do not rely on the results from the previous section. 

The definition of $w$ implies:
$$u_{x}=\beta w^{\beta-1}w_{x},\qquad
u_{xx}=\beta(\beta-1)w^{\beta-2}w_{x}^{2}+\beta w^{\beta-1}w_{xx},$$
$$u_{xxx}=\beta(\beta-1)(\beta-2)w^{\beta-3}w_{x}^{3} + 3\beta(\beta-1)w^{\beta-2}w_{x}w_{xx}+\beta w^{\beta-1}w_{xxx}.$$
Since $\beta \in (0,1)$, there are negative powers of $u$ in the above formulae. Nevertheless, so long as 
$u$ is bounded away from zero,  $u$ and $w$  live in the same $H^n$ spaces. Specifically:

\begin{lemma}\label{u-w}
Let $T>0$ be given.  Let $s\in\bf{N}$ such that  $s\geq 1$ be given.  Assume that $u$ is bounded uniformly away from zero.  Then, $u\in C([0,T];H^s)$ if and only if $w\in C([0,T];H^s).$\end{lemma}

This lemma shows that  in order to show existence of solutions for the initial value problem for $K(2,2),$ 
it is sufficient to show that the initial value problem for $w$ is solvable.

\noindent {\bf Proof:} We will only prove the backward direction of the implication; the proof of the forward implication is the same.
Since $u$ is bounded uniformly away from zero, $w$ is also bounded uniformly away from zero.
Let $\bar{c}>0$ be given such that for all $x\in X$ and for all $t\in[0,T],$ we have $w(x,t)>\bar{c}.$  

Since $w\in C([0,T];H^s)$ with $s\geq 1,$ by Sobolev imbedding we see $w\in C(X\times[0,T]),$ and therefore 
$u=w^{\beta}\in C(X\times[0,T]).$  This clearly implies $u\in C([0,T];L^2).$
Since $w(x,t)>\bar{c}>0$ for all $x$ and $t,$ we also have that $w^{\beta-1}\in C(X\times[0,T]).$  
Since $w\in C([0,T];H^s),$ we have $w_{x}\in C([0,T];H^{s-1})\subseteq C([0,T];L^2).$
Multiplying, we see that $w^{\beta-1}w_{x}\in C([0,T];L^2),$ so $u_{x}\in C([0,T];L^2),$ and therefore,
$u\in C([0,T];H^{1}).$

If $s\geq 2,$ then similarly, $w^{\beta-1}w_{xx}\in C([0,T];H^{s-2})\subseteq C([0,T];L^2).$  
Because $w(x,t)>\bar{c}>0$ for all $x$ and $t,$ we have that $w^{\beta-2}$ is continuous on $X\times[0,T].$
By the Sobolev algebra property, $w_{x}^{2}\in C([0,T];H^{s-1}),$ and therefore 
$w^{\beta-2}w_{x}^{2}\in C([0,T]; H^{s-1})\subseteq C([0,T]; L^2).$  From the above formula for $u_{xx},$ this implies
$u_{xx}\in C([0,T]; L^2),$ and therefore $u\in C([0,T]; H^{2}).$

We can continue in this manner for larger $s;$ this completes the proof.
\hfill $\blacksquare$

Now, we derive the evolution equation for $w.$  
Using (\ref{K22}) together with the definition of $w,$ we have
$$u_{t}=\beta w^{\beta-1} w_{t} = 2uu_{xxx} + 6u_{x}u_{xx} + 2uu_{x}.$$
We plug in the above formulas for $u$ and its first three derivatives, and solve for $w_{t},$ finding the following:
\begin{equation}\label{wEquation}
w_{t} = 2w^{\beta}w_{xxx} + (12\beta-6)w^{\beta-1}w_{x}w_{xx}+ (\beta-1)(8\beta-4)w^{\beta-2}w_{x}^{3}
+2w^{\beta}w_{x}.\end{equation}
We consider this equation with the initial data \begin{equation}\label{wData}w(x,0)=w_{0}(x).\end{equation}  
For $n\in{\bf N}$ with $n\geq 4,$ and for some positive
constant $a_0,$ we assume (similarly as in Assumption \ref{uniform}) the following:
\begin{equation}\label{wDataConditions}
w_{0}\in H^{n}\quad \mathrm{and}\quad w_{0}(x)\geq a_0,{\text{ for all }} x\in X.\end{equation}

The $K(2,2)$ equation has the following conserved quantities:
$$\int_{X} u^{3}(x,t)\ dx,\qquad \int_{X}u(x,t)\ dx.$$
The estimates we will be making can be viewed as generalizing the first of these.  In particular, we view
the first of these conserved quantities as being a conservation law for the $L^{2}$ norm of $u^{3/2};$
notice that in the case $n=0,$ we have $\beta=2/3,$ so $w=u^{3/2}.$  For larger values of $n,$ we find that
different powers of $u$ can be used to find a short-time estimate in $H^{n}.$

\subsection{Uniqueness and continuous dependence for $K(2,2)$}

In this section
we begin with a continuous dependence result
in $L^{2}.$  Uniqueness and continuous dependence in more regular spaces will be a corollary of this first result.
Since we are making estimates in $L^2$, we choose $\beta$ accordingly.
 Setting $n =0$ in \eqref{beta} gives
 $$\beta=2/3$$ and we take this value for $\beta$ throughout this subsection. 

\begin{theorem} \label{CDOIC-L2}
Let $w_{1}$ and $w_{2}$ be solutions of (\ref{wEquation}), with $w_{1}(x,0)=w_{0}(x)$ and 
$w_{2}(x,0)=\widetilde{w}_{0}(x),$ with both of these pieces of  
initial data satisfying (\ref{wDataConditions}).  Assume that there exists $T>0$ such that $w_{1}\in L^{\infty}([0,T]; H^{n})$
and $w_{2}\in L^{\infty}([0,T]; H^{n}).$  Assume also that for all $t\in[0,T],$ and for all $x\in X,$ we have
$w_{1}(x,t)>a_0/2$ and $w_{2}(x,t)>a_0/2.$
Then there exists $c>0$ such that for all $t\in [0,T],$
$$\|w_{1}(\cdot,t)-w_{2}(\cdot,t)\|_{L^2}^{2}\leq \left(\|w_{0}-\widetilde{w}_{0}\|_{L^2}^{2}\right)e^{ct}.$$
\end{theorem}

\noindent
{\bf Proof:}
Define $$E_{d}(t) = \frac{1}{2}\int_{X}(w_{1}(x,t)-w_{2}(x,t))^{2}\ dx.$$
Differentiating with respect to time, we find
$$\frac{dE_{d}}{dt} = \int_{X}(w_{1}-w_{2})(w_{1,t}-w_{2,t})\ dx.$$
We plug in from the evolution equation:
\begin{equation}\nonumber
\frac{dE_{d}}{dt} = A_{1} + A_{2} + A_{3} + A_{4},
\end{equation}
where 
\begin{equation}\nonumber
A_{1} = \int_{X} 2(w_{1}-w_{2})(w_{1}^{\beta}w_{1,xxx}-w_{2}^{\beta}w_{2,xxx})\ dx,\\
\end{equation}
\begin{equation}\nonumber
A_{2} = \int_{X}(12\beta-6)(w_{1}-w_{2})(w_{1}^{\beta-1}w_{1,x}w_{1,xx}-w_{2}^{\beta-1}w_{2,x}w_{2,xx})\ dx,
\end{equation}
\begin{equation}\nonumber
A_{3} = \int_{X}(\beta-1)(8\beta-4)(w_{1}-w_{2})(w_{1}^{\beta-2}w_{1,x}^{3}-w_{2}^{\beta-2}w_{2,x}^{3})\ dx,
\end{equation}
and 
\begin{equation}\nonumber
A_{4} = \int_{X}2(w_{1}-w_{2})(w_{1}^{\beta}w_{1,x}-w_{2}^{\beta}w_{2,x})\ dx.\end{equation}

We begin by adding and subtracting in $A_{1}:$
\begin{multline}\nonumber
A_{1}=\int_{X}2(w_{1}-w_{2})(w_{1}^{\beta}-w_{2}^{\beta})w_{1,xxx}\ dx 
\\+ \int_{X}2(w_{1}-w_{2})(w_{2}^{\beta})(w_{1,xxx}-w_{2,xxx})\ dx=:B_{1}+B_{2}.\end{multline}
On the domain of interest, the function $g(z)=z^{\beta}$ is Lipschitz continuous, so we have 
$$B_{1}\leq cE_{d}.$$
For $B_{2},$ we integrate by parts:
\begin{multline}\nonumber
B_{2} = -\int_{X}2(w_{1,x}-w_{2,x})(w_{2}^{\beta})(w_{1,xx}-w_{2,xx})\ dx 
\\-\int_{X}2\beta(w_{1}-w_{2})(w_{2}^{\beta-1}w_{2,x})(w_{1,xx}-w_{2,xx})\ dx.\end{multline}
We recognize that the first term on the right-hand side includes a perfect derivative, and we integrate by parts again:
\begin{multline}\nonumber
B_{2} = \int_{X}\beta(w_{2}^{\beta-1}w_{2,x})(w_{1,x}-w_{2,x})^{2}\ dx 
\\-\int_{X}2\beta(w_{1}-w_{2})(w_{2}^{\beta-1}w_{2,x})(w_{1,xx}-w_{2,xx})\ dx.\end{multline}
We integrate the second term on the right-hand side by parts:
\begin{multline}\nonumber
B_{2} = \int_{X}3\beta(w_{2}^{\beta-1}w_{2,x})(w_{1,x}-w_{2,x})^{2}\ dx 
\\+\int_{X}2\beta(w_{1}-w_{2})(w_{2}^{\beta-1}w_{2,x})_{x}(w_{1,x}-w_{2,x})\ dx.\end{multline}
For the second term on the right-hand side, we again notice a perfect derivative, and integrate by parts:
\begin{equation}\nonumber
B_{2} = \int_{X}3\beta(w_{2}^{\beta-1}w_{2,x})(w_{1,x}-w_{2,x})^{2}\ dx 
-\int_{X}\beta(w_{2}^{\beta-1}w_{2,x})_{xx}(w_{1}-w_{2})^{2}\ dx.\end{equation}

We now consider $A_{2}.$  To begin, we add and subtract:
\begin{multline}\nonumber
A_{2}=\int_{X}(12\beta-6)(w_{1}-w_{2})(w_{1}^{\beta-1}w_{1,x}-w_{2}^{\beta-1}w_{2,x})w_{1,xx}\ dx\\
+ \int_{X}(12\beta-6)(w_{1}-w_{2})(w_{2}^{\beta-1}w_{2,x})(w_{1,xx}-w_{2,xx})\ dx =: B_{3} + B_{4}.
\end{multline}
To estimate $B_{3},$ we add and subtract once more:
\begin{multline}\nonumber
B_{3}=\int_{X}(12\beta-6)(w_{1}-w_{2})(w_{1}^{\beta-1}-w_{2}^{\beta-1})w_{1,x}w_{1,xx}\ dx 
\\+ \int_{X}(12\beta-6)(w_{1}-w_{2})(w_{2}^{\beta-1})(w_{1,x}-w_{2,x})w_{1,xx}\ dx.\end{multline}
For the second term on the right-hand side, we recognize the presence of a perfect derivative, and we integrate by parts:
\begin{multline}\nonumber
B_{3}=\int_{X}(12\beta-6)(w_{1}-w_{2})(w_{1}^{\beta-1}-w_{2}^{\beta-1})w_{1,x}w_{1,xx}\ dx \\
- \int_{X}(6\beta-3)(w_{2}^{\beta-1}w_{1,xx})_{x}(w_{1}-w_{2})^{2}\ dx.\end{multline}
Using the uniform bounds, and the fact that the function $g(z)=z^{\beta-1}$ is Lipschitz on the relevant domain,
we find $$B_{3}\leq cE_{d}.$$
For $B_{4},$ we integrate by parts:
\begin{multline}\nonumber
B_{4} = -\int_{X}(12\beta-6)(w_{2}^{\beta-1}w_{2,x})(w_{1,x}-w_{2,x})^{2}\ dx
\\-\int_{X}(12\beta-6)(w_{1}-w_{2})(w_{2}^{\beta-1}w_{2,x})_{x}(w_{1,x}-w_{2,x})\ dx.
\end{multline}
For the second term on the right-hand side, we recognize a perfect derivative and we integrate by parts:
\begin{multline}\nonumber
B_{4} = -\int_{X}(12\beta-6)(w_{2}^{\beta-1}w_{2,x})(w_{1,x}-w_{2,x})^{2}\ dx
\\+\int_{X}(6\beta-3)(w_{2}^{\beta-1}w_{2,x})_{xx}(w_{1}-w_{2})^{2}\ dx.
\end{multline}
Now, since we have $\beta=2/3,$ and if we add $A_{1}+A_{2},$ there is an important cancellation, and we find
$$A_{1}+A_{2}=B_{1}+B_{2}+B_{3}+B_{4} \leq cE_{d}.$$

We can continue in the same manner (adding and subtracting, recognizing perfect derivatives, and integrating by parts),
and we will find $$A_{3}\leq cE_{d},\qquad A_{4}\leq cE_{d}.$$ 
Altogether, we have found $$\frac{dE_{d}}{dt}\leq cE_{d}.$$  This clearly implies
$$E_{d}(t)\leq E_{d}(0)e^{ct}.$$
This completes the proof of the lemma.\hfill$\blacksquare$

We can use this result to prove uniqueness, and we can also use it to get continuous dependence on the initial data in more
regular spaces.

\begin{cor}\label{K22 unique}
Let $w_{1}$ and $w_{2}$ be solutions of (\ref{wEquation}), with $w_{1}(x,0)=w_{0}(x)$ and 
$w_{2}(x,0)=\widetilde{w}_{0}(x),$ with the 
initial data satisfying (\ref{wDataConditions}).  Assume that there exists $T>0$ such that $w_{1}\in L^{\infty}([0,T]; H^{n})$
and $w_{2}\in L^{\infty}([0,T]; H^{n}).$  Assume also that for all $t\in[0,T],$ and for all $x\in X,$ we have
$w_{1}(x,t)>a_0/2$ and $w_{2}(x,t)>a_0/2.$  Let $n'$ be given such that $0\leq n'<n.$
Then, there exists a constant $C>0$ such that
for all $t\in[0,T],$ $$\|w_{1}(\cdot,t)-w_{2}(\cdot,t)\|_{H^{n'}}\leq C \|w_{0}-\widetilde{w}_{0}\|_{L^2}^{1-n'/n}.$$

\end{cor}

\noindent{\bf Proof:}  This is a straightforward combination of Theorem \ref{CDOIC-L2} and 
the Sobolev interpolation inequality used in \eqref{interp}.  
This completes the proof.
\hfill $\blacksquare$

\begin{cor}
Let $w_{1}$ and $w_{2}$ be solutions of (\ref{wEquation}), with $w_{1}(x,0)=w_{0}(x)$ and 
$w_{2}(x,0)=\widetilde{w}_{0}(x),$ with the 
initial data satisfying (\ref{wDataConditions}).  Assume that there exists $T>0$ such that $w_{1}\in L^{\infty}([0,T]; H^{n})$
and $w_{2}\in L^{\infty}([0,T]; H^{n}).$  Assume also that for all $t\in[0,T],$ and for all $x\in X,$ we have
$w_{1}(x,t)>a_0/2$ and $w_{2}(x,t)>a_0/2.$
If $w_{0}=\widetilde{w}_{0},$ then for all $t\in[0,T],$ $w_{1}(\cdot,t)=w_{2}(\cdot,t).$
\end{cor}

\noindent{\bf Proof:} This follows immediately from Theorem \ref{CDOIC-L2}.
\hfill $\blacksquare$

\

\noindent
\begin{remark} In light of Lemma \ref{u-w}, the uniqueness and continuous dependence theorems we have proved
for $w$ also give uniqueness and continuous dependence results for solutions of the original $K(2,2)$ equation,
(\ref{K22}).\end{remark}

\subsection{Existence of solutions for $K(2,2)$}

We now need to work in the regularity space of the initial data, $H^n$. Recall
$n\in{\bf N},$ and $n\geq 4.$  Take $\beta$ as in \eqref{beta}, that is $\beta=\displaystyle\frac{6}{2n+9}.$

We introduce a mollifed version of \eqref{wEquation}.
Specifically take  $\mathcal{J}_\ep$, for $\epsilon > 0$, to be the same
mollifier as was described in Theorem \ref{mollifier stuff}.
 Our mollified equation is
\begin{multline}\label{evolution}
{w^{\epsilon}_{t}} = \mathcal{J}_{\epsilon}\Bigg[
2(\mathcal{J}_{\epsilon}{w^{\epsilon}})^{\beta}(\mathcal{J}_{\epsilon}{w^{\epsilon}_{xxx}})
+(12\beta-6)(\mathcal{J}_{\epsilon}{w^{\epsilon}})^{\beta-1}(\mathcal{J}_{\epsilon}{w^{\epsilon}_x})
(\mathcal{J}_{\epsilon}{w^{\epsilon}_{xx}})\\
+ (\beta-1)(8\beta-4)(\mathcal{J}_{\epsilon}{w^{\epsilon}})^{\beta-2}
(\mathcal{J}_{\epsilon}{w^{\epsilon}_x})^{3}
+2(\mathcal{J}_{\epsilon}{w^{\epsilon}})^{\beta}(\mathcal{J}_{\epsilon}{w^{\epsilon}_x})\Bigg].\end{multline}
This evolution equation is taken with the initial condition
\begin{equation}\label{initialCondition}
{w^{\epsilon}}(x,0)=w_{0}(x),\end{equation}
where $w_{0}\in H^{n}$ and there exists $a_0>0$ such that for all $x\in X,$ 
$w_{0}(x)>a_0.$
We will seek solutions of this initial value problem which satisfy at positive times $t$ the condition
\begin{equation}\label{awayFromZero}
{w^{\epsilon}}(x,t)>\frac{a_0}{2}\ \text{for all} \ x\in X.\end{equation}

Given $w_{0}\in H^{n}$ which satisfies the condition $w_{0}(x)>a_0$ for all $x\in X,$ we define an open set
$\mathcal{O}\subset H^n$ such that $w_{0}\in\mathcal{O}.$
Specifically
$$
\mathcal{O}:=\left\{ f \in H^n : \|f\|_{H^n}< 2\|w_{0}\|_{H^n}\ \text{and, for all} \ x\in X,\ |f(x)|>\frac{a_0}{2} \right\}.
$$
The following existence result is a direct application of the Picard theorem, and we note it without further proof.

\begin{lemma}\label{PicardApplication}
For any $\epsilon>0,$ there exists $T_{\epsilon}>0$ such that there exists 
${w^{\epsilon}}(x,t)\in C^{1}([-T_{\epsilon},T_{\epsilon}]; \mathcal{O})$ such that ${w^{\epsilon}}$ 
is a solution of the initial value problem (\ref{evolution}), (\ref{initialCondition}).
\end{lemma}
\begin{remark}
To be clear, for all $t\in [-T_{\epsilon}, T_{\epsilon}],$ since ${w^{\epsilon}}(\cdot,t)\in\mathcal{O},$
the solution ${w^{\epsilon}}$ satisfies the condition (\ref{awayFromZero}).
\end{remark}
\

The time of existence in Lemma \ref{PicardApplication} depends badly on $\epsilon.$  In order to remove this dependence,
we will perform an $H^{n}$ energy estimate for this equation; this will allow us to use the continuation theorem for 
autonomous differential equations.  Our energy is the following:
$$E(t)=\frac{1}{2}\int_{X}\left({w^{\epsilon}}\right)^{2}\ dx + \frac{1}{2}\int_{X}(\partial_{x}^{n}{w^{\epsilon}})^{2}\ dx := E_{0} + E_{n}.$$
Since we are estimating ${w^{\epsilon}}$ in $H^{n},$ it will be helpful to apply $n$ spatial derivatives to (\ref{evolution}):
\begin{multline}\label{evolutionWithDerivatives}
\partial_{x}^{n}{w^{\epsilon}_{t}} = \mathcal{J}_{\epsilon}\Bigg[2(\mathcal{J}_{\epsilon}{w^{\epsilon}})^{\beta}
\partial_{x}^{n+3}\mathcal{J}_{\epsilon}{w^{\epsilon}}
+ (12\beta+2\beta n-6)(\mathcal{J}_{\epsilon}{w^{\epsilon}})^{\beta-1}
(\mathcal{J}_{\epsilon}{w^{\epsilon}_x})\partial_{x}^{n+2}\mathcal{J}_{\epsilon}{w^{\epsilon}}
\\+F_{1}[\mathcal{J}_{\epsilon}{w^{\epsilon}}]\partial_{x}^{n+1}\mathcal{J}_{\epsilon}{w^{\epsilon}} 
+F_{2}[\mathcal{J}_{\epsilon}{w^{\epsilon}}]\Bigg].\end{multline}
Formulas and estimates for $F_{1}$ and $F_{2}$ are the subject of the next subsection.

\subsubsection{The Lower-Order Terms}

The formula for $F_{1}$ is
\begin{multline}
F_{1}[\mathcal{J}_{\epsilon}{w^{\epsilon}}]=2{n \choose 2}\left((\mathcal{J}_{\epsilon}{w^{\epsilon}})^{\beta}\right)_{xx}
+ (12\beta-6)(\mathcal{J}_{\epsilon}{w^{\epsilon}})^{\beta-1}\mathcal{J}_{\epsilon}{w^{\epsilon}_{xx}}\\
+(12\beta-6)n\partial_{x}\left((\mathcal{J}_{\epsilon}{w^{\epsilon}})^{\beta-1}\mathcal{J}_{\epsilon}{w^{\epsilon}_x}\right)
+3(\beta-1)(8\beta-4)(\mathcal{J}_{\epsilon}{w^{\epsilon}})^{\beta-2}(\mathcal{J}_{\epsilon}{w^{\epsilon}_x})^{2}
+2(\mathcal{J}_{\epsilon}{w^{\epsilon}})^{\beta}.
\end{multline}
Upon expanding the various derivatives in the above formula, we find
$$F_{1}=d_{1}(\mathcal{J}_{\epsilon}{w^{\epsilon}})^{\beta-2}(\mathcal{J}_{\epsilon}{w^{\epsilon}_x})^{2}
+d_{2}(\mathcal{J}_{\epsilon}{w^{\epsilon}})^{\beta-1}(\mathcal{J}_{\epsilon}{w^{\epsilon}})_{xx}
+2(\mathcal{J}_{\epsilon}{w^{\epsilon}})^{\beta},$$
for some constants $d_{1}$ and $d_{2}.$

Recall that we are considering ${w^{\epsilon}}$ such that for all $x$ and $t,$ ${w^{\epsilon}}(x,t)>0.$  Since $\mathcal{J}_{\epsilon}$ is an averaging
operator, for any $z\in X$ and for any $t$ we have
$\displaystyle\inf_{x\in X} {w^{\epsilon}}(x,t)\leq \mathcal{J}_{\epsilon}{w^{\epsilon}}(z,t),$ and therefore 
$$\left\|\frac{1}{\mathcal{J}_{\epsilon}{w^{\epsilon}}(z)}\right\|_{L^\infty} \leq \left\|\frac{1}{{w^{\epsilon}}}\right\|_{L^{\infty}}.$$
Since $\beta<1,$ we have $\beta-1<0,$ and $\beta-2<0,$ and so on.  We then estimate, for instance, 
$$|\mathcal{J}_{\epsilon}{w^{\epsilon}}|^{\beta-1}\leq \left\|\frac{1}{{w^{\epsilon}}}\right\|_{L^{\infty}}^{1-\beta}.$$
We estimate other negative powers of $\mathcal{J}_{\epsilon}{w^{\epsilon}}$ in the same way.

We therefore have the following bounds for $F_{1}:$
$$|F_{1}|_{L^{\infty}}\leq \gamma_{1}\left(\left\|\frac{1}{{w^{\epsilon}}}\right\|_{L^{\infty}}\right)
\left(\|{w^{\epsilon}}\|_{3}^{2}+\|{w^{\epsilon}}\|_{3}^{\beta}\right)\leq \gamma_{1}\left(\left\|\frac{1}{{w^{\epsilon}}}\right\|_{L^{\infty}}\right)
\left(E + E^{\beta/2}\right),$$
\begin{equation}\label{F1xEstimate}
|\partial_{x}F_{1}|_{L^{\infty}} \leq \gamma_{2}\left(\left\|\frac{1}{{w^{\epsilon}}}\right\|_{L^{\infty}}\right)
\left(\|{w^{\epsilon}}\|_{4}^{3}+\|{w^{\epsilon}}\|_{4}\right)
\leq \gamma_{2}\left(\left\|\frac{1}{{w^{\epsilon}}}\right\|_{L^{\infty}}\right)\left(E^{3/2}+E^{1/2}\right),\end{equation}
for some smooth functions $\gamma_{1}$ and $\gamma_{2}.$  Note: we only actually will need the second of these two bounds.

For $F_{2},$ we clearly have the definition
\begin{multline}
F_{2}[\mathcal{J}_{\epsilon}{w^{\epsilon}}] =
\partial_{x}^{n}\Bigg[
2(\mathcal{J}_{\epsilon}{w^{\epsilon}})^{\beta}(\mathcal{J}_{\epsilon}{w^{\epsilon}_{xxx}})
+(12\beta-6)(\mathcal{J}_{\epsilon}{w^{\epsilon}})^{\beta-1}(\mathcal{J}_{\epsilon}{w^{\epsilon}_x})
(\mathcal{J}_{\epsilon}{w^{\epsilon}_{xx}})\\
+ (\beta-1)(8\beta-4)(\mathcal{J}_{\epsilon}{w^{\epsilon}})^{\beta-2}
(\mathcal{J}_{\epsilon}{w^{\epsilon}_x})^{3}
+2(\mathcal{J}_{\epsilon}{w^{\epsilon}})^{\beta}(\mathcal{J}_{\epsilon}{w^{\epsilon}_x})\Bigg]\\
-
\Bigg[2(\mathcal{J}_{\epsilon}{w^{\epsilon}})^{\beta}
\partial_{x}^{n+3}\mathcal{J}_{\epsilon}{w^{\epsilon}}
+ (12\beta+2\beta n-6)(\mathcal{J}_{\epsilon}{w^{\epsilon}})^{\beta-1}
(\mathcal{J}_{\epsilon}{w^{\epsilon}_x})\partial_{x}^{n+2}\mathcal{J}_{\epsilon}{w^{\epsilon}}
\\+F_{1}[\mathcal{J}_{\epsilon}{w^{\epsilon}}]\partial_{x}^{n+1}\mathcal{J}_{\epsilon}{w^{\epsilon}} \Bigg].
\end{multline}
For ease of estimating, it is helpful to rewrite this using the product rule:
\begin{multline}\label{F2Sums}
F_{2}[\mathcal{J}_{\epsilon}{w^{\epsilon}}] =
2\sum_{k=3}^{n}{n \choose k}\left(\partial_{x}^{k}\left((\mathcal{J}_{\epsilon}{w^{\epsilon}})^{\beta}\right)\right)
(\partial_{x}^{n-k+3}\mathcal{J}_{\epsilon}{w^{\epsilon}})\\
+(12\beta-6)\sum_{k=2}^{n-1}{n \choose k}
\left(\partial_{x}^{k}\left((\mathcal{J}_{\epsilon}w^{\epsilon})^{\beta-1}\mathcal{J}_{\epsilon}w^{\epsilon}_{x}\right)\right)
\left(\partial_{x}^{n-k+2}\mathcal{J}_{\epsilon}w^{\epsilon}\right)
\\
+(12\beta-6)\sum_{k=1}^{n}{n \choose k}(\mathcal{J}_{\epsilon}w^{\epsilon}_{xx})
\left(\partial_{x}^{k}(\mathcal{J}_{\epsilon}w^{\epsilon})^{\beta-1}\right)
\left(\partial_{x}^{n-k+1}\mathcal{J}_{\epsilon}w^{\epsilon}\right)
\\+3(\beta-1)(8\beta-4)\sum_{k=1}^{n-1}{n-1 \choose k}
\left(\partial_{x}^{k}\left((\mathcal{J}_{\epsilon}w^{\epsilon})^{\beta-2}
(\mathcal{J}_{\epsilon}w^{\epsilon}_{x})^{2}\right)\right)(\partial_{k}^{n+1-k}w^{\varepsilon})
\\+(\beta-1)(8\beta-4)\sum_{k=0}^{n-1}{n-1 \choose k}
\left(\partial_{x}^{k+1}\left(\mathcal{J}_{\epsilon}w^{\epsilon}\right)^{\beta-2}\right)
\left(\partial_{x}^{n-1-k}\left(\mathcal{J}_{\epsilon}w^{\epsilon}_{x}\right)^{3}\right)
\\+2\sum_{k=1}^{n}{n \choose k}
\left(\partial_{x}^{k}\left((\mathcal{J}_{\epsilon}{w^{\epsilon}})^{\beta}\right)\right)
\left(\partial_{x}^{n-k+1}\mathcal{J}_{\epsilon}{w^{\epsilon}}\right).
\end{multline}

We now give an estimate for $F_{2}$ in $L^2.$  For this estimate, note that the highest number of derivatives on
${w^{\epsilon}}$ that appear anywhere in (\ref{F2Sums}) is $n.$  So, $F_{2}$ will be bounded in terms of $E,$ but we 
should be careful as to which powers of $E$ can be used for the bound, and we need to be careful here because we
estimate positive powers and negative powers differently.  Notice that in the fifth summation on the right-hand side of
(\ref{F2Sums}), the $k=n-1$ term includes a factor 
$\partial_{x}^{n}\left((\mathcal{J}_{\epsilon}{w^{\epsilon}})^{\beta-2}\right)$.
In this factor, if we apply the derivatives by using the chain rule and the product rule, then we find one term which has
as a factor $(\mathcal{J}_{\epsilon}{w^{\epsilon}})^{\beta-2-n}.$ This is the most negative power encountered of all the 
terms in (\ref{F2Sums}).  Similarly, we want to identify the least negative power on the right-hand side of (\ref{F2Sums});
this comes, for instance, 
from the $k=1$ term in the sixth summation, and is $(\mathcal{J}_{\epsilon}{w^{\epsilon}})^{\beta-1}.$
From these considerations, we find the following estimate:
\begin{equation}\label{F2Estimate}
\|F_{2}\|_{L^2}\leq \gamma_{3}\left(\left\|\frac{1}{{w^{\epsilon}}}\right\|_{L^{\infty}}\right)
\left(\|{w^{\epsilon}}\|_{H^n}^{n+3}+\|{w^{\epsilon}}\|_{H^n}^{2}\right)
\leq \gamma_{3}\left(\left\|\frac{1}{{w^{\epsilon}}}\right\|_{L^{\infty}}\right)\left(E^{(n+3)/2}+E\right),\end{equation}
for some smooth function $\gamma_{3}.$

\subsubsection{The Energy Estimate}

The goal of this section is to prove the following lemma:
\begin{lemma}
Let $\epsilon>0$ be given.
Assume there exists $\bar{T}>0$ such that there exists ${w^{\epsilon}}\in C([0,\bar{T}]; H^{n})$ which satisfies
(\ref{evolution}) and (\ref{awayFromZero}).  Then there exist positive constants $C_{1}$ and $C_{2}$ which depend
only on $a_0$ and $n$ such that for all $t\in[0,\bar{T}],$ we have
$$E(t)\leq -\frac{\ln\left(\exp\{-C_{2}E(0)\}-C_{1}C_{2}t\right)}{C_{2}}.$$
\end{lemma}

\noindent{\bf Proof:}
The method of proof is to first prove an estimate for the growth of the energy.  This provides a differential inequality, and
the solution of this inequality gives the estimate we are trying to establish.

Clearly, the growth of $E_{0}$ is bounded by $E:$
$$\frac{dE_{0}}{dt} \leq \gamma_{4}\left(\left\|\frac{1}{{w^{\epsilon}}}\right\|_{L^{\infty}}\right)\left(E^{1+\beta/2}+E^{2}\right),$$
for some smooth $\gamma_{4}.$

We consider now the growth of $E_{n}:$
\begin{multline}\label{Es}
\frac{dE_{n}}{dt} = F + \int_{X} 
2(\mathcal{J}_{\epsilon}{w^{\epsilon}})^{\beta}
(\partial_{x}^{n}\mathcal{J}_{\epsilon} {w^{\epsilon}})
(\partial_{x}^{n+3}\mathcal{J}_{\epsilon}{w^{\epsilon}})\\
+(12\beta+2\beta n-6)(\mathcal{J}_{\epsilon}{w^{\epsilon}})^{\beta-1}
(\mathcal{J}_{\epsilon}{w^{\epsilon}_x})
(\partial_{x}^{n}\mathcal{J}_{\epsilon}{w^{\epsilon}})
(\partial_{x}^{n+2}\mathcal{J}_{\epsilon}{w^{\epsilon}})\ dx,
\end{multline}
where we have introduced the term
$$F=\int_{X} \left(F_{1}[\mathcal{J}_{\epsilon}{w^{\epsilon}}]\partial_{x}^{n+1}\mathcal{J}_{\epsilon}{w^{\epsilon}} 
+F_{2}[\mathcal{J}_{\epsilon}{w^{\epsilon}}]\right)(\partial_{x}^{n}\mathcal{J}_{\epsilon}{w^{\epsilon}})\ dx.$$
We can bound $F$ in terms of the energy easily; first, this requires one integration by parts:
$$F=\int_{X} -\frac{1}{2}(F_{1}[\mathcal{J}_{\epsilon}{w^{\epsilon}}])_{x}
(\partial_{x}^{n}\mathcal{J}_{\epsilon}{w^{\epsilon}})^{2} 
+\left(F_{2}[\mathcal{J}_{\epsilon}{w^{\epsilon}}]\right)(\partial_{x}^{n}\mathcal{J}_{\epsilon}{w^{\epsilon}})\ dx.$$
By the above bounds (\ref{F1xEstimate}) and (\ref{F2Estimate}) for the $F_{i},$ and since $n\geq 4,$ we have
$$|F|\leq \gamma_{5}\left(\left\|\frac{1}{{w^{\epsilon}}}\right\|_{L^{\infty}}\right)\left(E^{(n+4)/2}+E^{3/2}\right),$$
for some smooth $\gamma_{5}.$

We integrate by parts in the first term on the right-hand side of (\ref{Es}):
\begin{multline}\label{Es2}
\frac{dE_{n}}{dt} = F+ \int_{X} \Bigg[-2(\mathcal{J}_{\epsilon}{w^{\epsilon}})^{\beta}
(\partial_{x}^{n+1}\mathcal{J}_{\epsilon}{w^{\epsilon}})
(\partial_{x}^{n+2}\mathcal{J}_{\epsilon}{w^{\epsilon}}) \\
+(10\beta+2\beta n-6)(\mathcal{J}_{\epsilon}{w^{\epsilon}})^{\beta-1}
(\mathcal{J}_{\epsilon}{w^{\epsilon}_x})
(\partial_{x}^{n}\mathcal{J}_{\epsilon}{w^{\epsilon}})
(\partial_{x}^{n+2}\mathcal{J}_{\epsilon}{w^{\epsilon}})\Bigg]\ dx.\end{multline}
We recognize that the first term in the integrand on the right-hand side of (\ref{Es2}) includes a perfect derivative, 
and we integrate it by parts:
\begin{multline}\label{Es3}
\frac{dE_{n}}{dt} = F +  \int_{X} \Bigg[\beta (\mathcal{J}_{\epsilon}{w^{\epsilon}})^{\beta-1}
(\mathcal{J}_{\epsilon}{w^{\epsilon}_x})(\partial_{x}^{n+1}\mathcal{J}_{\epsilon}{w^{\epsilon}})^{2}  \\
+(10\beta+2\beta n-6)(\mathcal{J}_{\epsilon}{w^{\epsilon}})^{\beta-1}
(\mathcal{J}_{\epsilon}{w^{\epsilon}_x})(\partial_{x}^{n}\mathcal{J}_{\epsilon}{w^{\epsilon}})
(\partial_{x}^{n+2}\mathcal{J}_{\epsilon}{w^{\epsilon}})\Bigg]\ dx.\end{multline}
Now, we integrate the second term in the integrand on the right-hand side of (\ref{Es3}) by parts:
\begin{multline}\label{Es4}
\frac{dE_{n}}{dt} = F + \int_{X}\Bigg[ 
(-9\beta-2\beta n+6) (\mathcal{J}_{\epsilon}{w^{\epsilon}})^{\beta-1}
(\mathcal{J}_{\epsilon}{w^{\epsilon}_x})(\partial_{x}^{n+1}\mathcal{J}_{\epsilon}{w^{\epsilon}})^{2}
\\-(10\beta+2\beta n-6)
((\mathcal{J}_{\epsilon}{w^{\epsilon}})^{\beta-1}(\mathcal{J}_{\epsilon}{w^{\epsilon}_x}))_{x}
(\partial_{x}^{n}\mathcal{J}_{\epsilon}{w^{\epsilon}})
(\partial_{x}^{n+1}\mathcal{J}_{\epsilon}{w^{\epsilon}})\Bigg]\ dx.\end{multline}
Since $\beta=\displaystyle\frac{6}{2n+9},$ the first term in the integrand on the right-hand side of (\ref{Es4}) is zero
(this is the reason for selecting $\beta$ this way):
\begin{equation}\label{Es5}
\frac{dE_{n}}{dt} = F 
-\int_{X}
(10\beta+2\beta n-6)((\mathcal{J}_{\epsilon}{w^{\epsilon}})^{\beta-1}(\mathcal{J}_{\epsilon}{w^{\epsilon}_x}))_{x}
(\partial_{x}^{n}\mathcal{J}_{\epsilon}{w^{\epsilon}})(\partial_{x}^{n+1}\mathcal{J}_{\epsilon}{w^{\epsilon}})
\ dx.\end{equation}
We integrate by parts in the right-hand side of (\ref{Es5}) once more:
\begin{equation}\label{Es6}
\frac{dE_{n}}{dt} = 
F+\int_{X}(5\beta+\beta n-3)((\mathcal{J}_{\epsilon}{w^{\epsilon}})^{\beta-1}(\mathcal{J}_{\epsilon}{w^{\epsilon}_x}))_{xx}
(\partial_{x}^{n}\mathcal{J}_{\epsilon}{w^{\epsilon}})^{2}\ dx.\end{equation}
Finally, we have 
\begin{equation}\label{almostMainEstimate}
\frac{dE}{dt}\leq \gamma_{6}\left(\left\|\frac{1}{{w^{\epsilon}}}\right\|_{L^{\infty}}\right)
\left(E^{(n+4)/2}+E^{3/2}\right),\end{equation}
for some smooth $\gamma_{6}.$

Since condition (\ref{awayFromZero}) is satisfied, we can rewrite (\ref{almostMainEstimate}) as
\begin{equation}\label{mainEstimate}
\frac{dE}{dt}\leq \bar{c}\left(E^{(n+4)/2}+E^{3/2}\right).\end{equation}
There exist positive constants $C_{1}$ and $C_{2}$ such that 
$$\bar{c}\left(E^{(n+4)/2}+E^{3/2}\right)\leq C_{1}\exp\{C_{2}E\},$$
and these constants depend only on $a_0$ and $s.$  So, we have
$$\frac{dE}{dt}\leq C_{1}\exp\{C_{2}E\}.$$
This differential inequality can be solved, with the result that
$$E(t)\leq -\frac{\ln\left(\exp(-C_{2}E(0))-C_{1}C_{2}t\right)}{C_{2}}.$$
This completes the proof.\hfill$\blacksquare$

\subsubsection{Passage to the limit}

We now introduce another open set, $\widetilde{\mathcal{O}} \subset H^n$ which is defined as:
$$
\widetilde{\mathcal{O}}:= \left\{ f \in H^n : \|f\|_{H^n}< \frac{3}{2}\|w_{0}\|_{H^n}\ \text{and, for all } x\in X,\ |f(x)|>\frac{3a_0}{4}.\right\}
$$
Clearly, $\widetilde{\mathcal{O}}\subseteq\mathcal{O}.$

\begin{lemma}\label{uniformTimeInterval}
Let $w_{0}\in H^{n}$ satisfy (\ref{awayFromZero}).
There exists $T^{*}>0$ such that for all $\epsilon\in(0,1),$ the solution ${w^{\epsilon}}$ of the initial value problem
(\ref{evolution}), (\ref{initialCondition}) satisfies ${w^{\epsilon}}\in C([0,T^{*}]; \widetilde{\mathcal{O}}).$
\end{lemma}

\noindent{\bf Proof:}
Assume there does not exist such a $T^{*}.$  This implies that there exists a sequence $\epsilon_{j}\rightarrow 0$ such that
the solutions $w_{\epsilon_{j}}$ leave the set $\widetilde{\mathcal{O}}$ arbitrarily quickly as $j\rightarrow \infty.$  
For each $j,$ let
$T_{j}$ be the infimum of the set of times at which $w_{\epsilon_{j}}$ is not in $\widetilde{\mathcal{O}}.$ 
From Lemma \ref{PicardApplication}, we know $T_{j}>0$ 
for all $j,$ and by choice of the sequence, we have $\displaystyle\lim_{j\rightarrow\infty}T_{j}=0.$
So, for each $j,$ we have $w_{\epsilon_{j}}\in C([0,T_{j}); \widetilde{\mathcal{O}}),$
and $w_{\epsilon_{j}}(\cdot,T_{j})\notin\widetilde{\mathcal{O}}.$ Note, however, that we do have 
$w_{\epsilon_{j}}\in C([0,T_{j}];\mathcal{O}).$  In particular, we must have either 
$$\inf_{x\in X}w_{\epsilon_{j}}(x,T_{j})=\frac{3a_0}{4},\qquad \mathrm{or}
\qquad E(T_{j})= \frac{9}{4}E(0).$$
But this is impossible, from the above uniform estimate.  In particular, if we let
 $$\tilde{T} = \frac{\exp(-C_{2}E(0))-\exp(-2C_{2}E(0))}{C_{1}C_{2}},$$
then we see that if $T_{j}\in[0,\tilde{T}],$ then we have $$E(T_{j})\leq 2E(0).$$
So, we cannot have $E(T_{j})=\frac{9}{4}E(0)$ with $T_{j}\rightarrow 0.$

From the evolution equation for ${w^{\epsilon}},$ we see that there exists $C>0$ such that for all $\epsilon\in(0,1),$
we have $|{w^{\epsilon}_{t}}|\leq C$ as long as $w^{\epsilon}\in\mathcal{O}.$  Therefore,
if the infimum of ${w^{\epsilon}}$ was originally at least $a_0,$ it could not become equal to $\frac{3a_0}{4}$ arbitrarily quickly.
This completes the proof. \hfill $\blacksquare$

\begin{theorem}\label{K22 exist}
Let $w_{0}\in H^{n}$ satisfy (\ref{awayFromZero}).  Then, there exists $T>0$ such that the initial value problem 
(\ref{evolution}), (\ref{initialCondition}) has a solution $w\in C([0,T]; \mathcal{O}).$
\end{theorem}

\noindent
{\bf Proof:} We let $T$ equal the value $T^{*}$ from Lemma \ref{uniformTimeInterval}.  For $\epsilon\in(0,1),$ 
the functions ${w^{\epsilon}}$ are all in the space $C([0,T]; \mathcal{O}),$ which implies that
for all $\epsilon\in(0,1),$ for all $t\in[0,T],$ we have the uniform bound $\|{w^{\epsilon}}(\cdot, t)\|_{H^n}\leq 2\|w_{0}\|_{H^n}.$
Since $s\geq 4,$ this implies that there exists $C>0$ such that for all $\epsilon\in(0,1),$ for all $t\in[0,T],$ 
for all $x\in X,$ we have
$$|{w^{\epsilon}}(x,t)|\leq C,\qquad |{w^{\epsilon}_x}(x,t)|\leq C,\qquad |{w^{\epsilon}_{t}}(x,t)|\leq C.$$
This means that the family ${w^{\epsilon}}$ is uniformly bounded and equicontinuous, on the domain $X\times[0,T].$
The Arzela-Ascoli theorem gives a uniform limit of the sequence ${w^{\epsilon}}$ 
(taking a subsequence, which we do not relabel)
in $C(X\times [0,T]).$  Together with the uniform bound and an elementary interpolation inequality,
this implies that the limit, $w,$ is in $C([0,T]; H^{n'}),$ for any $0\leq n'<n.$  

Furthermore, pointwise in time, ${w^{\epsilon}}$ is uniformly bounded in $H^{n},$ so there is a weak limit 
of a subsequence (since the unit ball
of a Hilbert space is weakly compact) in $H^{n}.$  By uniqueness of limits, this implies that for all $t,$ 
we have $w(\cdot, t)\in H^{n}.$
We now want to show that $w\in C([0,T];H^{n}),$ and that $w$ solves the non-mollified equation.

To show that $w$ solves the non-mollified evolution equation, we integrate the mollified evolution equation in time:
\begin{multline}\label{integratedInTime}
{w^{\epsilon}}(x,t) = w_{0}(x) + \int_{0}^{t}2\mathcal{J}_{\epsilon}\left(
\left(\left(\mathcal{J}_{\epsilon}\left({w^{\epsilon}}\right)^{\beta}\right)(x,s)\right)
\mathcal{J}_{\epsilon}{w^{\epsilon}_{xxx}}(x,s)\right)\ ds \\
+\int_{0}^{t}(12\beta-6)\mathcal{J}_{\epsilon}\left(
\left(\left(\mathcal{J}_{\epsilon}{w^{\epsilon}}\right)^{\beta-1}(x,s)\right)
\left(\mathcal{J}_{\epsilon}{w^{\epsilon}_x}(x,s)\right)
\mathcal{J}_{\epsilon}{w^{\epsilon}_{xx}}(x,s)\right)\ ds\\
+\int_{0}^{t}(\beta-1)(8\beta-4)\mathcal{J}_{\epsilon}\left(
\left(\left(\mathcal{J}_{\epsilon}{w^{\epsilon}}\right)^{\beta-2}(x,s)\right)
\left(\mathcal{J}_{\epsilon}{w^{\epsilon}_x}(x,s)\right)^{3}
\right)\ ds \\
+ \int_{0}^{t}2\mathcal{J}_{\epsilon}\left(\left(\left(\mathcal{J}_{\epsilon}{w^{\epsilon}}\right)^{\beta}(x,s)\right)
\mathcal{J}_{\epsilon}{w^{\epsilon}_x}(x,s)\right)\ ds,
\end{multline}
where $t\in[0,T].$
Since ${w^{\epsilon}}$ converges uniformly to $w$ in $C([0,T];H^{n'})$ for any $0\leq n'<n,$ and since $n\geq 4,$
we can pass to the limit in (\ref{integratedInTime}):
\begin{multline}\nonumber
w(x,t) = w_{0}(x) +\\
 \int_{0}^{t}2w^{\beta}(x,s)w_{xxx}(x,s) + (12\beta-6)w^{\beta-1}(x,s)w_{x}(x,s)w_{xx}(x,s)\ ds\\
+\int_{0}^{t}(\beta-1)(8\beta-4)w^{\beta-2}(x,s)w_{x}^{3}(x,s) + 2w^{\beta}(x,s)w_{x}(x,s)\ ds.
\end{multline}
This clearly implies that $w$ is a solution of the initial value problem.

All that remains is to demonstrate the highest regularity.  We will begin by showing that $\|w\|_{H^n}$ is continuous as a function 
of time.  We start by showing that the norm is right-continuous at $t=0.$  By Fatou's Lemma,
we have $$\|w(\cdot,0)\|_{H^n}^{2}\leq \liminf_{t\rightarrow 0^{+}}\|w(\cdot,t)\|_{H^n}^{2}.$$
From the energy inequality, however, we have 
$$\limsup_{t\rightarrow 0^{+}} \|w(\cdot,t)\|_{H^n}^{2}\leq \|w(\cdot,0)\|_{H^n}^{2}.$$
Therefore, the norm is right-continuous at $t=0.$

Given $t^{*}\in(0,T),$ we consider $t^{*}$ as a new initial time.  We may solve the initial value problem starting from
time $t^{*},$ with initial data $w(\cdot, t^{*}).$   By the previous argument, we see that this solution is right-continuous at 
time $t=t^{*}.$  By the uniqueness theorem, we know that this solution starting from time $t^{*}$
and the previous solution, $w,$ are the same.  Therefore, $\|w\|_{H^n}$ is right-continuous at all $t\in[0,T).$
Furthermore, all of the analysis we have performed works with time reversed: there is no feature of our estimates which
requires time to move forward.  Therefore, all of this analysis could be reversed, finding that $\|w\|_{H^n}$ is left-continuous
at all $t\in(0,T].$

To complete the proof, we use the fact that weak convergence together with convergence of the norm implies strong convergence.
We have already shown convergence of the norm, so as the final step in our proof, we show that $w$ is weakly continuous in time.
Let $\psi\in H^{-n}$ be given.  Given $t^{*}\in[0,T],$ we want to show 
$$\langle \psi, w(\cdot, t) - w(\cdot, t^{*})\rangle_{L^2} \rightarrow 0,$$
as $t\rightarrow t^{*}.$  Let $\delta>0$ be given.  Let $K>0$ be such that for all $t\in[0,T],$ $\|w(\cdot,t)\|_{H^n}\leq K.$
Let $n'$ be given such that $0\leq n'<n.$  Then, of course, $-n<-n',$ and $H^{-n'}$ is dense in $H^{-n}.$  Therefore we can find
$\psi_{\delta}$ such that $$\|\psi_{\delta}-\psi\|_{-n}\leq \frac{\delta}{3(1+K)}.$$
This implies the following bound:
$$\left|\langle \psi-\psi_{\delta}, w(\cdot,t)-w(\cdot,t^{*})\rangle_{L^2}\right| \leq \left(\frac{\delta}{3(1+K)}\right)
\left(2K\right)\leq \frac{2\delta}{3}.$$
Now, since $w\in C([0,T]; H^{n'}),$ we can take $t$ close enough to $t^{*}$ such that 
$$\|w(\cdot,t)-w(\cdot,t^{*})\|_{H^{n'}}\leq \frac{\delta}{3(1+\|\psi_{\delta}\|_{-n'})}.$$
This implies the following bound:
$$\left|\langle \psi_{\delta}, w(\cdot,t)-w(\cdot,t^{*})\rangle_{L^2}\right| \leq \|\psi_{\delta}\|_{-n'}
\left(\frac{\delta}{3(1+\|\psi_{\delta}\|_{-n'})}\right) \leq \frac{\delta}{3}.$$
Putting this together, we find that we can take $t$ sufficiently close to $t^{*}$ to get 
$$\left|\langle \psi, w(\cdot, t) - w(\cdot, t^{*})\rangle_{L^2}\right| \leq \delta.$$
This implies that $w$ is weakly continuous in time.  This completes the proof.
\hfill$\blacksquare$

We sum up our results for the $K(2,2)$ equation in the following corollary:
\begin{cor}\label{K22exist}
Let $n$ be an integer such that $n\geq 4.$  Let $X$ be a periodic interval.  Let $a_0>0.$
Let $u_{0}\in H^{n}$ be such that $u_{0}(x)>a_0$ for all $x\in X.$
Then there exists $T>0$ and $u\in C([0,T]; H^{n})$ such that for all $x\in X$ and for all $t\in[0,T],$ $u$ 
satisfies $u(x,t)>a_0/2,$ and $u$ is the unique solution of
 the initial value problem (\ref{K22}) with $u(x,0)=u_{0}(x).$  
 
 Let $n'$ be given such that $0\leq n' < n,$ and let $\epsilon>0$ be given. There exists $\delta>0$ such that
 for any $v_{0}\in H^{n}$ such that $v_{0}(x)> a_0$ for all $x\in X$
 and such that $\|v_{0}-u_{0}\|_{H^n}\leq \delta,$
 the solution, $v,$ of the initial value problem (\ref{K22}) with $v(x,0)=v_{0}(x)$ 
satisfies $v\in C([0,T]; H^{n}),$ and
$$\sup_{t\in[0,T]}\|u(\cdot,t)-v(\cdot,t)\|_{H^{n'}} < \epsilon.$$
\end{cor}

\section{Other examples of quasilinear equations}

In this section, we give additional examples of quasilinear equations for which our method gives a short-time
well-posedness result.  Like (\ref{K22}), these are also equations in which the dispersive effect appears nonlinearly.
The first additional example is the Harry Dym equation, and the other examples arise in numerical analysis of finite
difference schemes.

\subsection{Well-posedness of the Harry Dym Equation}

The Harry Dym equation is 
\begin{equation}\label{HD}
u_{t} = u^{3}u_{xxx}.\end{equation}
This is a completely integrable equation \cite{kruskal} which has been shown to have
applications in interfacial fluid dynamics \cite{tanveer}.
We consider this with the initial condition
\begin{equation}\label{HDIC}u(x,0)=u_{0}(x).\end{equation}
As in the case of the $K(2,2)$ equation, we can show that the initial value problem
for the Harry Dym equation is well-posed for $H^{4}$ initial data which is bounded away from
zero.  The argument is completely analogous to the $K(2,2)$ case, so we only present the 
essential step, which is the $H^{n}$ energy estimate.  When providing all the details of the
well-posedness proof, as we have done for the $K(2,2)$ equation, the energy estimate must
in fact be carried out for a mollified version of the equation; in this section, we proceed 
informally, and present the energy estimate without first introducing mollifiers.  Certainly, however,
all of the details of the present argument can be carried out carefully as in the preceding section.
We remark that, as in the case of the $K(2,2)$ equation, the average modified diffusion, $\bar{\delta},$ that we defined 
before is identically zero; in the case of the Harry Dym equation, this is because there is no $u_{xx}$ term
present in (\ref{HD}).

Two conserved quantities for the Harry Dym equation are the following:
$$\int_{X}\frac{1}{u(x,t)}\ dx,\qquad \int_{X}\frac{1}{u^{2}(x,t)}\ dx.$$
One view of the energy estimate which we are about to perform 
is that we are generalizing the first of these.
We can take the view that the first of these conserved quantities tells us that the $L^{2}$ norm of
$u^{-1/2}$ is conserved by the evolution.

Let $n$ be an integer, such that $n\geq 4.$  We make the choice $\beta=\displaystyle
\frac{1}{n-\frac{1}{2}},$ following the same construction
as used for \eqref{K22} above. Note that, as before, this implies $\beta\in(0,1);$
furthermore, we also have $3\beta\in(0,1).$
As before, we define $w=u^{\frac{1}{\beta}},$ so that $u=w^{\beta}.$
From the previous section, we have the following relevant formulas:
$$u_{t}=\beta w^{\beta-1}w_{t},$$
$$u_{xxx}=\beta(\beta-1)(\beta-2)w^{\beta-3}w_{x}^{3}
+3\beta(\beta-1)w^{\beta-2}w_{x}w_{xx} + \beta w^{\beta-1}w_{xxx}.$$
Combining these, we find the evolution equation for $w:$
\begin{equation}\label{dym-w}
w_{t} = w^{3\beta}w_{xxx}
+3(\beta-1)w^{3\beta-1}w_{x}w_{xx}
+(\beta-1)(\beta-2)w^{3\beta-2}w_{x}^{3}.
\end{equation}

We apply the operator $\partial_{x}^{n}$ to (\ref{dym-w}):
\begin{equation}
\partial_{x}^{n}w_{t} = w^{3\beta}\partial_{x}^{n+3}w
+(3\beta + 3\beta n - 3)w^{3\beta -1}w_{x}\partial_{x}^{n+2}w
+F_{1}[w]\partial_{x}^{n+1}w+F_{2}[w],
\end{equation}
where
\begin{multline}\nonumber
F_{1}[w] = {n \choose 2}\partial_{x}^{2}(w^{3\beta})
+ 3(\beta - 1)n(w^{3\beta-1}w_{x})_{x}
\\+ 3(\beta - 1)w^{3\beta-1}w_{xx}
+3(\beta-1)(\beta-2)w^{3\beta-2}w_{x}^{2},
\end{multline}
and
\begin{multline}\nonumber
F_{2}[w] = \sum_{k=3}^{n}{n \choose k}(\partial_{x}^{k}w^{3\beta})(\partial_{x}^{n-k+3}w)
+\sum_{\ell=1}^{n}{n\choose\ell}(\partial_{x}^{\ell}w^{3\beta-1})(\partial_{x}^{n-\ell+1}w)w_{xx}
\\
+3(\beta-1)\sum_{k=2}^{n-1}\sum_{\ell=1}^{k}{n\choose k}{k\choose\ell}(\partial_{x}^{\ell}w^{3\beta-1})
(\partial_{x}^{k-\ell+1}w)(\partial_{x}^{n-k+2}w)
\\
+ \left(\sum_{k=1}^{s-1}(\partial_{x}^{k}(3w^{3\beta-2}w_{x}^{2}))(\partial_{x}^{n+1-k}w)\right)
+\partial_{x}^{n-1}((\partial_{x}w^{3\beta-2})w_{x}^{3}).
\end{multline}
Since $F_{1}[w]$ includes at most second derivatives of $w,$ we can bound the $L^{\infty}$ norm 
of $F_{1}[w]$ in terms of $\|w\|_{H^n}.$  Since $F_{2}[w]$ includes at most $n$-many derivatives of
$w,$ we can bound the $L^{2}$ norm of $F_{2}[w]$ in terms of $\|w\|_{H^n}.$

We proceed to the energy estimate.  We define the energy, $E(t),$ to be 
$$E(t)= \frac{1}{2}\int_{X} w^{2} + (\partial_{x}^{n}w)^{2}\ dx.$$
Then, differentiating, we find the following:
$$\frac{dE}{dt} = \int_{X} ww_{t} + (\partial_{x}^{n}w)(\partial_{x}^{n}w_{t})\ dx.$$
Clearly, since $n\geq 4,$ we can bound $\displaystyle\int_{X}ww_{t}\ dx$ in terms of $E.$
Also, the contributions from $F_{1}$ and $F_{2}$ are bounded in terms of $E.$

The remaining piece that we must concern ourselves with is 
$$\int_{X}(\partial_{x}^{n}w)w^{3\beta}\partial_{x}^{n+3}w
+(3\beta + 3\beta n - 3)(\partial_{x}^{n}w)w^{3\beta -1}w_{x}\partial_{x}^{n+2}w\ dx.$$
We integrate the first term by parts:
$$\int_{X}-(\partial_{x}^{n+1}w)w^{3\beta}\partial_{x}^{n+2}w 
+(3\beta n-3)(\partial_{x}^{n}w)w_{x}\partial_{x}^{n+2}w\ dx.$$
We recognize a perfect derivative in the first term, and we integrate both terms by parts:
$$\int_{X}\left(\frac{3\beta}{2}-3\beta n + 3\right)w^{3\beta-1}w_{x}(\partial_{x}^{n+1}w)^{2}
-(3\beta n -3)(\partial_{x}^{n}w)w_{xx}\partial_{x}^{n+1}w\ dx.$$
Since $\beta=\displaystyle\frac{1}{n-\frac{1}{2}},$ the first of these terms is identically zero.
Performing one more integration by parts, we see that the final term is bounded in terms of the energy.

\begin{theorem}\label{HD exist}
Let $n$ be an integer such that $n\geq 4.$  Let $X$ be a periodic interval.  Let $a_0>0.$
Let $u_{0}\in H^{n}$ be such that $u_{0}(x)>a_0$ for all $x\in X.$
Then there exists $T>0$ and $u\in C([0,T]; H^{n})$ such that for all $x\in X$ and for all $t\in[0,T],$ $u$ 
satisfies $u(x,t)>a_0/2,$ and $u$ is the unique solution of
 the initial value problem (\ref{HD}), (\ref{HDIC}).  
 
 Let $n'$ be given such that $0\leq n' < n,$ and let $\epsilon>0$ be given. There exists $\delta>0$ such that
 for any $v_{0}\in H^{n}$ such that $v_{0}(x)> a_0$ for all $x\in X$
 and such that $\|v_{0}-u_{0}\|_{H^n}\leq \delta,$
 the solution, $v,$ of the initial value problem (\ref{HD}) with $v(x,0)=v_{0}(x)$ 
satisfies $v\in C([0,T]; H^{n}),$ and
$$\sup_{t\in[0,T]}\|u(\cdot,t)-v(\cdot,t)\|_{H^{n'}} < \epsilon.$$
\end{theorem}

\subsection{Effective equations for finite difference schemes}

In the paper \cite{goodmanLax}, Goodman and Lax found that a certain finite difference scheme for the 
KdV equation has as an modified equation 
\begin{equation}\label{g-l}
u_{t}+uu_{x}+\frac{1}{6}\Delta^{2}uu_{xxx}=0,\end{equation}
where $\Delta$ is the constant spatial step.
Similarly, in studying a finite-difference scheme for a nonlocal dispersive equation, Zumbrun found the following
modified equation \cite{zumbrun}:
\begin{equation}\label{z}
u_{t}+(u^{2})_{x}=-c_{2}a^{2}(uu_{xx})_{x},
\end{equation}
where $c_{2}$ and $a$ are constants.

Slight variations of the arguments of the previous sections, which we will not repeat, imply that both of these are well-posed in $H^{4},$ for
initial data bounded away from zero.  The important point is that 
as for the $K(2,2)$ equation and the Harry Dym equation, the average modified diffusion for (\ref{g-l}) and for (\ref{z}) is identically zero.

\section{Remarks on singularity formation}

We mention some further results related to singularity formation; these results hold for all of the quasilinear equations we have studied. 
These singularity formation results stem from theorems in \cite{ambroseWright1} on positivity preservation for equations with the form
$$u_{t}=uF[u]+u_{x}G[x],$$
where $F$ and $G$ are operators which can include $u$ and derivatives of $u.$  (We note that the paper \cite{deFrutos} also contains
a positivity preservation theorem for the $K(2,2)$ equation.)  The following theorem is a version of Theorem 6 of 
\cite{ambroseWright1}, which is proved using ideas from \cite{constantinEscher}:

\begin{theorem}Let $T>0$ and let $u\in C([0,T]; H^{4}(X))$ be a solution of (\ref{K22}),  (\ref{HD}), (\ref{g-l}), or  (\ref{z}).   
For each $t\in[0,T],$ let $m(t)=\inf_{x\in X} u(x,t).$  Then, 
for all $t\in[0,T],$ we have $\mathrm{sgn}(m(t))=\mathrm{sgn}(m(0)).$\end{theorem}

The following corollary is an immediate consequence of this:

\begin{cor} Let $n\in{\bf{N}}$ with $n\geq 4.$  Let $a_{0}>0.$  Let $u_{0}\in H^{n}(X)$ be such that 
$u_{0}(x)>a_{0},$ for all $x\in X.$  Let $u$ be the solution of the initial value problem for (\ref{K22}), (\ref{HD}),
(\ref{g-l}), or (\ref{z}), with $u(x,0)=u_{0}(x).$  Assume there exists $x_{*}\in X$ and $t_{*}>0$ such that
$u(x_{*},t_{*})=0.$  Then there exists $t_{**}\in(0,t_{*}]$ such that $$\lim_{t\rightarrow t_{**}^{-}}\|u(\cdot, t)\|_{4} = +\infty.$$
\end{cor}

\bibliography{harry-dym}{}

\begin{thebibliography}{10}

\bibitem{Akhunov}
T.~Akhunov.
\newblock {\em Local well posedness of dispersive systems in 1D}.
\newblock PhD thesis, University of Chicago, 2011.

\bibitem{ASWY}
D.~M. Ambrose, G.~Simpson, J.~D. Wright, and D.~G. Yang.
\newblock Ill-posedness of degenerate dispersive equations.
\newblock Preprint, 2011.

\bibitem{ambroseWright1}
David~M. Ambrose and J.~Douglas Wright.
\newblock Preservation of support and positivity for solutions of degenerate
  evolution equations.
\newblock {\em Nonlinearity}, 23(3):607--620, 2010.

\bibitem{chertock}
Alina Chertock and Doron Levy.
\newblock Particle methods for dispersive equations.
\newblock {\em J. Comput. Phys.}, 171(2):708--730, 2001.

\bibitem{constantinEscher}
Adrian Constantin and Joachim Escher.
\newblock Wave breaking for nonlinear nonlocal shallow water equations.
\newblock {\em Acta Math.}, 181(2):229--243, 1998.

\bibitem{Craig-Kappeler-Strauss}
W.~Craig, T.~Kappeler, and W.~Strauss.
\newblock Gain of regularity for equations of {K}d{V} type.
\newblock {\em Ann. Inst. H. Poincar\'e Anal. Non Lin\'eaire}, 9(2):147--186,
  1992.

\bibitem{Craig-Goodman}
Walter Craig and Jonathan Goodman.
\newblock Linear dispersive equations of {A}iry type.
\newblock {\em J. Differential Equations}, 87(1):38--61, 1990.

\bibitem{deFrutos}
J.~de~Frutos, M.~{\'A}. L{\'o}pez~Marcos, and J.~M. Sanz-Serna.
\newblock A finite difference scheme for the {$K(2,2)$} compacton equation.
\newblock {\em J. Comput. Phys.}, 120(2):248--252, 1995.

\bibitem{evans}
Lawrence~C. Evans.
\newblock {\em Partial differential equations}, volume~19 of {\em Graduate
  Studies in Mathematics}.
\newblock American Mathematical Society, Providence, RI, second edition, 2010.

\bibitem{goodmanLax}
Jonathan Goodman and Peter~D. Lax.
\newblock On dispersive difference schemes. {I}.
\newblock {\em Comm. Pure Appl. Math.}, 41(5):591--613, 1988.

\bibitem{kruskal}
Martin Kruskal.
\newblock Nonlinear wave equations.
\newblock In {\em Dynamical systems, theory and applications ({R}encontres,
  {B}attelle {R}es. {I}nst., {S}eattle, {W}ash., 1974)}, pages 310--354.
  Lecture Notes in Phys., Vol. 38. Springer, Berlin, 1975.

\bibitem{levyShuYan}
Doron Levy, Chi-Wang Shu, and Jue Yan.
\newblock Local discontinuous {G}alerkin methods for nonlinear dispersive
  equations.
\newblock {\em J. Comput. Phys.}, 196(2):751--772, 2004.

\bibitem{liOlverRosenau}
Yi~A. Li, Peter~J. Olver, and Philip Rosenau.
\newblock Non-analytic solutions of nonlinear wave models.
\newblock In {\em Nonlinear theory of generalized functions ({V}ienna, 1997)},
  volume 401 of {\em Chapman \& Hall/CRC Res. Notes Math.}, pages 129--145.
  Chapman \& Hall/CRC, Boca Raton, FL, 1999.

\bibitem{majdaBertozzi}
Andrew~J. Majda and Andrea~L. Bertozzi.
\newblock {\em Vorticity and incompressible flow}, volume~27 of {\em Cambridge
  Texts in Applied Mathematics}.
\newblock Cambridge University Press, Cambridge, 2002.

\bibitem{Rosenau-Hyman}
Philip Rosenau and J.M. Hyman.
\newblock Compactons: solitons with finite wavelength.
\newblock {\em Phys. Rev. Lett.}, (70):564--567, 1993.

\bibitem{rus-pade}
Francisco Rus and Francisco~R. Villatoro.
\newblock Pad\'e numerical method for the {R}osenau-{H}yman compacton equation.
\newblock {\em Math. Comput. Simulation}, 76(1-3):188--192, 2007.

\bibitem{rus-selfsimilar}
Francisco Rus and Francisco~R. Villatoro.
\newblock Self-similar radiation from numerical {R}osenau-{H}yman compactons.
\newblock {\em J. Comput. Phys.}, 227(1):440--454, 2007.

\bibitem{tanveer}
S.~Tanveer.
\newblock Evolution of {H}ele-{S}haw interface for small surface tension.
\newblock {\em Philos. Trans. Roy. Soc. London Ser. A}, 343(1668):155--204,
  1993.

\bibitem{zumbrun}
Kevin Zumbrun.
\newblock On a nonlocal dispersive equation modeling particle suspensions.
\newblock {\em Quart. Appl. Math.}, 57(3):573--600, 1999.

\end{thebibliography}
\bibliographystyle{plain}
\end{document}